\documentclass[11pt]{article}
\usepackage{amsfonts,amssymb, graphicx,a4,bm,bbm}

\newcommand{\zeq}{\setcounter{equation}{0}}

\newtheorem{theo}{Theorem}[section]
\newtheorem{lem}[theo]{Lemma}
\newtheorem{prop}[theo]{Proposition}
\newtheorem{defi}[theo]{Definition}

\let\a=\alpha \let\b=\beta \let\ch=\chi \let\d=\delta \let\e=\varepsilon
 \let\g=\gamma     \let\k=\kappa 
\let\m=\mu  \let\o=\omega     
\let\r=\rho \let\s=\sigma \let\t=\tau 
 \let\x=\xi 
\let\D=\Delta \let\F=\Phi \let\G=\Gamma  
\let\O=\Omega    

\let\\=\noindent
\def\ee{\end{equation}}
\def\be{\begin{equation}}
\def\vv{\vskip.2cm}
\def\0{\emptyset}
\def\T{{\mathcal{T}}}
\def\PP{{\mathcal P}}
\def\ef{\mathfrak{f}}

\title{A local lemma via entropy compression}
\author{
Rog\'erio G. Alves$^{1}$,
Aldo Procacci$^2$ and
Remy Sanchis$^2$\\
\\
\small{$^1$ Departamento de Matem\'atica UFOP}
\small{ 35400-000 - Ouro Preto - MG
Brazil}
\\
\small{$^2$ Departamento de Matem\'atica UFMG}
\small{ 30161-970 - Belo Horizonte - MG
Brazil}}

\begin{document}

\maketitle

\begin{abstract}
In the framework of the probabilistic method in combinatorics, we revisit the entropy compression method  clarifying  the setting
in which it can be applied and  providing a theorem yielding a general constructive criterion. We finally  elucidate,  through topical examples,  the effectiveness of the
entropy-compression criterion in comparison with the  Lovasz Local Lemma criterion and, in particular,  with the improved criterion based on cluster expansion.

\end{abstract}

\vskip.2cm
{\footnotesize

}
\vskip.5cm

\section{Introduction}
The Lov\'asz Local Lemma (LLL), originally formulated by Erd\H{o}s and
Lov\'{a}sz \cite{EL}, is a powerful tool in the framework of the
probabilistic method in combinatorics to prove the existence of
some combinatorial objects with certain desirable properties (such as
a proper coloring of the edges of a graph). According  the LLL,  the existence of the combinatorial object
under analysis is guaranteed if a family of (bad) events
 in some probability space do not occur. The LLL then provides a criterion providing an upper bound  on the probabilities of these bad events
which ensures that, with non-zero probability, none of them occur. The  popularity of this lemma is
due to the fact that it can  be implemented
for a wide class of problems  in combinatorics,  in such a way that its
sufficient condition, once some few parameters  have been suitably tuned,
can be easily checked.

\def\fa{{\mathfrak{E}}}\def\ev{{\mathfrak{e}}}
\def\ev{{\mathfrak{e}}}
 \def\fa{{\mathfrak{F}}} \def\GG{{\cal G}}\def\prob{{\rm P}} \def\N{{\mathbb{N}}}
\subsection{ Lov\'asz Local Lemma: the general  version}
%In order to enunciate  the Lov\'asz local Lemma  we need to introduce the concept of dependency graph of a family of events.
In order to enunciate the most general statement   of the LLL, which includes also the lopsided version  (see e.g.  \cite{SS} and references therein), we need to give some basic notations.
Let $G=(V,E)$ be a (simple undirected) graph with vertex set $V$ and edge set $E$.
Two vertices $v,v'\in V$ are adjacent if  $\{v,v'\}\in E$. Given $v\in V$, let $\G^*_G(v)$ be
the set of vertices of $G$ adjacent to $v$ (i.e., the neighborhood of $v$).
We denote   $\G_G(v)= \G^*_G(v)\cup\{v\}$.

\\Given a family of events $\fa$ in some probability space $(\O,{\rm P})$ and given $\ev\in \fa$, we denote by  $\overline \ev$  the complement
event of $\ev\in \fa$, so that   $\bigcap_{\ev\in \fa}\overline \ev$ is  the event that none of the events
in the family $\fa$  occur.

\\Henceforth,  the product over the empty set  is equal to one,
if $n\in \N$ we set $[n]\equiv\{1,2,\dots, n\}$ and
$[n]_0\equiv\{0,1,2,\dots, n\}$ and, if $X$
is a finite set, then $|X|$ denotes the number of its elements.

\begin{theo}[Lov\'asz local Lemma]\label{lov}
Let $\fa$ be a finite collection of events in a probability space $(\O,\prob)$.
Let $\GG$ be a graph with vertex set  $\fa$.  Let  $\bm{\mu}=\{\mu_\ev\}_{\ev\in \fa}$ be a collection of real numbers
 in $(0,+\infty)$. If, for each event  $\ev\in \fa$ and for each ${U}\subset \fa\setminus \G_\GG(\ev)$,
\begin{equation}
\label{cond}
\prob\left(\ev\Big|\bigcap_{\ev'\in U}\bar{\ev'}\right)\;\leq\;  {\mu_\ev \over  \Phi_\ev(\bm \m,\GG)}
\end{equation}
with
\be\label{Phi}
\Phi_\ev(\bm \m,\GG)= \prod_{\ev'\in \Gamma_\GG(\ev)} (1+\mu_{\ev'})
\ee
then
\be\label{thesis}
\prob\Big(\bigcap_{\ev\in \fa}\overline \ev\Big)\,> \,0
\ee
i.e. the probability that none of the events in the family $\cal \fa$  occur is strictly positive.
\end{theo}
The sufficent condition given by (\ref{cond}) is the LLL criterion.
\vskip.2cm

\\{\bf Remark}. Theorem \ref{lov} is generally known in the literature as the ``Lopsided" LLL.
When the  graph $\GG$ has  edge set
such that,
each event $\ev\in \fa$   is independent of
the $\s$-algebra generated by the collection of events $\fa \setminus \Gamma_\GG(\ev)$ so that, for any   $\ev\in \fa$ and any  ${U}\subset \fa\setminus \G_\GG(\ev)$
we have that $\prob\left(\ev\Big|\bigcap_{\ev'\in U}\bar{\ev'}\right)= \prob(\ev)$, we say that $\GG$ is a {\it dependency graph} for the family $\fa$. Theorem \ref{lov} in which $\GG$ is a dependency graph for the family $\fa$ is the most common version of the LLL and it is the one used
in the vast majority of the applications. An even more simpler version of the LLL is the so-called ``symmetric" case.
This occurs if it is possible to find  $p>0$ such that  $\prob\big(\ev\big|\cap_{\ev'\in U}\bar{\ev'}\big)\le p$ for all events in $\fa$
and any  ${U}\subset \fa\setminus \G_\GG(\ev)$. In this case it is natural to choose $\m_\ev=\mu$ for mall $\ev\in \fa$ so that condition (\ref{cond}) becomes
\begin{equation}
\label{condsym}
p\;\leq\; {\mu \over (1+\m)^{d+1}}
\end{equation}
where $d=\max_{\ev\in \fa} |\G^*_\GG(\ev)|$.  Optimizing in $\mu>0$ and using that $\Big(1+{1\over d}\Big)^d<e$, (\ref{condsym}) is usually written
as
\be
\label{condsym2}
ep(d+1)\;\leq\; 1
\end{equation}
%We will use  the LLL criterion in the simple form (\ref{condsym2}) in Section \ref{PropB}.

\subsection{The improved version of the LLL via cluster expansion}
In recent times a remarkable  connection between the LLL and statistical mechanics  has been
pointed out by Scott and Sokal \cite{SS}.
More precisely Scott and Sokal showed that the thesis (\ref{thesis})
of the Lov\'asz local lemma holds  {\it if and only if}
the probabilities $\{\prob(\ev)\}_{\ev\in \mathfrak{E}}$
are inside the complex polydisc where the partition function
of the hard-core self-repulsive lattice gas on
$G$ is nonvanishing and that this necessary and sufficient condition is equivalent to the Shearer ``optimal"
criterion \cite{Sh}. Furthermore, they show that the
LLL criterion (\ref{cond})) providing a {\it sufficient} condition for  the thesis  (\ref{thesis}) to hold,
coincides  with  the so-called Dobrushin criterion, which is a sufficient condition for the non vanishing
of the aforementioned partition function.

\\Later on, Fern\'andez and Procacci  \cite{FP}, working on the abstract polymer gas model (of which the
hard-core, self-repulsive lattice gas on
a graph $G$ is a special case)
obtained a sensible improvement of the Dobrushin criterion
via cluster expansion methods. Such a result has then been used by Bissacot et al.  \cite{BFPS}  to obtain an
improved version of the LLL. To enunciate this
improved version of the LLL we recall that, given a graph $G=(V, E)$, a set
$Y\subset V$  is {\it independent} if no edge  $\{v,v'\}\in E$ is such that $v\in Y$ and $v'\in Y$.

%\\The new version of the LLL can then be stated as follows.

\begin{theo}[Cluster Expansion local lemma (CELL)]\label{bis}
With the same hypothesis of Theorem \ref{lov},  if, for each event  $\ev\in \fa$ and for each ${U}\subset \fa\setminus \G_\GG(\ev)$
\begin{equation}
\label{cond1n}
\prob\left(\ev\Big|\bigcap_{\ev'\in U}\bar{\ev'}\right)\;\leq\;  {\mu_\ev \over  \Xi_\ev(\bm \m,\GG)}
\end{equation}
with
\be\label{Xi}
\Xi_\ev(\bm \m,\GG)=\sum_{Y\subseteq\Gamma_\GG(\ev)\atop Y~{\rm independent}} \prod_{\ev'\in Y}\mu_{\ev'}
\ee
then
$$
{\prob}\Big(\bigcap_{\ev\in \fa}\bar \ev\Big)\,> \,0
$$
\end{theo}
\\{\bf Remark}.
The improvement on Theorem \ref{lov} is immediately recognized by noting that
$$
\Phi_\ev(\bm \m,\GG)= \prod_{\ev'\in \G_\GG(\ev)} (1+\mu_{\ev'})= \sum\limits_{Y\subseteq  \Gamma_\GG(\ev)}
\prod_{\ev'\in Y}\mu_{\ev'} \ge \sum\limits_{Y\subseteq  \Gamma_\GG(\ev)\atop Y\ {\rm indep\ in}\ \GG}
\prod_{\ev'\in Y}\mu_{\ev'}=\Xi_\ev(\bm \m,\GG)
$$ \def\A{{{\mathcal O}}}

\\The improved  sufficient condition given by (\ref{cond1n}) is nowadays known as the {\it cluster expansion  criterion} (CE criterion).
As remarked in \cite{FP},  observing that
the event $\ev$ is connected to any $\ev'\in \G^*_\GG(\ev)$
the function  $\Xi_\ev(\bm \m,\GG)$  can be  rewritten as
\be\label{Xi2}
\Xi_\ev(\bm \m,\GG)=\mu_\ev+\sum_{Y\subseteq\Gamma^*_\GG(\ev)\atop Y~{\rm independent}} \prod_{\ev'\in Y}\mu_{\ev'}
\ee
%We will use  (\ref{Xi2}) in Section \ref{PropB}.

\subsection{The Moser-Tardos algorithmic version of the local lemma}
The LLL provides a sufficient condition for the probability that
none of the undesirable events in some (unspecified) probability $(\O,\prob)$ space occur to be strictly positive  and this
implies the existence of at least one outcome in $(\O,\prob)$
which realizes the occurrence of the ``good" event.  However, since  $(\O,\prob)$  is left completely unspecified,   Theorem \ref{lov},
as well as Theorem \ref{bis}  says nothing about how to find
such a configuration. The efforts to devise an  algorithmic version of the LLL by trying to say something more precise
about the probability space  $(\O,\prob)$   go back to the work of Beck
 \cite{B},  Alon \cite{A} and others   (see e.g. reference list in  \cite{Mo2}) and
 culminated with the  breakthrough work by Moser and Tardos \cite{Mo,MT}, who  presented a fully  algorithmic version
 of LLL by assuming that the   probability space $(\O,\prob)$  is generated by a finite collection of independent random variables. This simple assumption
 covers nearly  all known applications of LLL.
\def\obj{{\rm obj}}

\\The scheme   introduced by    Moser and Tardos  (called nowadays the {\it  variable setting}) is defined as follows.
  Given a finite set $\mathbf{S}$, for each  $x\in \mathbf{S}$ we have
a random variable   $\psi_x$ taking values in some set $\Psi_x$ according with some given distribution.
The collection $\{\psi_x\}_{x\in \mathbf{S}}$  is supposed to form a family of  {\it mutually independent} random variables.
Therefore, for any non empty $U\subset \mathbf{S}$, the subfamily  $\{\psi_x\}_{x\in U}$  generates a product probability space
  $\O_U$  with product probability measure $\prob_U$.
  Moser and Tardos set  $(\O,\prob)=(\O_\mathbf{S}, \prob_\mathbf{S})$.

\\Observe that  for any $U\subset \mathbf{S}$, a configuration (or elementary event) $\o=\prod_{x\in U}\psi_x$ in $\O_U$ is an element
of  $\Psi_U=\prod_{x\in U}\Psi_x$.
Given $\o\in \Psi_\mathbf{S}$ and $U\subset  \mathbf{S} $,
 we let $\o|_{U}$ be the restriction of $\o$ to $U$, i.e. if $\o=\prod_{x\in \mathbf{S}}\psi_x$ then  $\o|_{U}=\prod_{x\in U}\psi_x$.

  \\Let now  ${{\mathcal O}}$
 be   a  family of subsets of $\mathbf{S}$ (called hereafter ``objects"). Given $A\in \A$,
an event $\ev$ in $\O_A$ is a subset of $\Psi_U=\prod_{x\in U}\Psi_x$.
For each $A\in \A$, we will refer to an event in $\O_A$  as  {\it an event associated to $A$} and write ${\rm obj}(\ev)=A$.
As usual, $\bar \ev$ denotes the complementary event
of $\ev$ . Clearly, since $\O_\mathbf{S}$
is a product space, each  event  $\ev$ in $\O_A$  can also be viewed as
the  event in $\O_\mathbf{S}$   formed by those $\o\in \Psi_\mathbf{S}$ such that the restriction $\o|_A$ is in $ \ev$
so that
$\prob_A(\ev)=\prob_\mathbf{S}(\ev)$.

\\Moser and Tardos considered families $\fa$ of (bad) events such that,
for each element of  $\ev\in \fa$,  $\obj(\ev)=A$,  for some $A\in {\mathcal O}$.

\vskip.2cm
\\{\bf Remark}. Note
 that by construction an $\ev\in \fa$  depends  only on the family of  random variables
$\{\psi_x\}_{x\in \obj(\ev)}$ and   any two events $\ev, \ev'$ of $\fa$ such that $\obj(\ev)\cap \obj(\ev')= \emptyset$
are necessarily independent.\def\E{{\cal E}}
This implies that
the graph $\GG$  with vertex set $\fa$ and edge-set  constituted by the pairs
$\{\ev,\ev'\}$ such that $\obj(\ev)\cap \obj(\ev')\neq \emptyset$ is  a natural dependency graph for the family  $\fa$.

\vskip.2cm

\\In this setting  Moser and Tardos defined the following algorithm. \def\A{{{\mathcal O}}}
\vskip.2cm
{\small \\{\bf MT-Algorithm}.
\vskip.1cm
\\- Step 0. Sample a random  $\o\in \Psi_\mathbf{S}$ and set $\o_0=\o$ as the configuration at step $0$.

\\- Step $i$ (for $i\ge 1$). Let $\o_{i-1}\in \Psi_\mathbf{S}$ be the configuration at step $i-1$. If $\o_{i-1}$ is such that  some event
of the family  $\fa$ occurs,
select one  of them (at random or according to some deterministic rule), say $\ev$,
   and  choose an  evaluation (resampling) $\o'\in \Psi_{\obj(\ev)}$. Set
$\o_{i}$ as the configuration of $\Psi_\mathbf{S}$ such that $\o_{i}|_{\mathbf{S}\setminus \obj(\ev)}=  {\o_{i-1}}|_{\mathbf{S}\setminus  \obj(\ev)}$ and
${\o_{i}}|_{ \obj(\ev)}= \o'$
}
\vskip.2cm
\\Moser and Tardos proved  \cite{MT}
that if condition (\ref{cond}) of  Theorem \ref{lov} holds,  MT-algorithm  terminates rapidly finding
a configuration in  $\Psi_\mathbf{S}$ such that none of the bad events of the family $\fa$ occurs. Later Pegden \cite{Pe}
improved Moser-Tardos result  replacing
condition (\ref{cond}) with condition (\ref{cond1n}). Pegden's extension  of the Moser-Tardos algorithmic local lemma  is thus the algorithmic version of
Theorem \ref{bis} and can be stated as follows.

\begin{theo}[{\bf Algorithmic  CELL}]\label{MosTar} Given  a finite set $\mathbf{S}$  and its associated family
of mutually independent random variables $\psi_\mathbf{S}$,
let ${\mathcal O}$ be a collection  of subsets of $\mathbf{S}$ and  let $\fa$ a family of bad events associated to $\A$ with natural dependency graph
$\GG$.
Let  $\bm{\mu}=\{\mu_\ev\}_{\ev\in \fa}$ be  real numbers
 in $(0,+\infty)$. If, for each $\ev\in \fa$, the cluster expansion criterion (\ref{cond1n}) holds,
then
$$
\bigcap_{\ev\in \fa}\bar \ev\neq\0
$$
and
 MT-algorithm finds $\o\in\bigcap_{\ev\in \fa}\bar \ev$ in an  expected total number of  steps less than or equal to
$\sum_{\ev\in \fa_\A}\m_\ev$.
\end{theo}
{\bf Remark}. In the Moser-Tardos setting
above described the function $\Xi_\ev(\bm \m,\GG)$ defined in (\ref{Xi}) admits
a somehow natural  upper bound.
Just observe that, for  any $y\in {\mathbf{S}}$, the set   $\fa(y)=\{\ev\in \fa:~y\in \obj(\ev)\}$
is a clique of the natural dependency graph $\GG$ of the family $\fa$: all events in $\fa(y)$ contain  $y$ an thus any pair $\{\ev,\ev'\}\subset \fa(y)$ is connected by an edge of $\GG$.
Thus, the neighbor $\G_\GG(\ev)$ of any event  $\ev\in \fa$ is  the (in general not disjoint) union of cliques $\{\fa(y)\}_{y\in \obj(\ev)}$, i.e., we can write
$\G_\GG(\ev)= \cup_{y\in \obj(\ev)}\fa(y)$.
Therefore
\be\label{clique}
\Xi_\ev(\bm \m,\GG)\le  \Xi^{\rm clique}_\ev(\bm \m,\GG)\equiv\prod_{y\in \obj(\ev)}\Big[1+ \sum_{\ev'\in \fa(y)}\m_{\ev'}\Big]
\ee
A sligtly better clique estimate is obtained by using the expression (\ref{Xi2}) for $ \Xi_\ev(\bm \m,\GG)$. Namely,
\be\label{cliqueb}
\Xi_\ev(\bm \m,\GG)\le \widetilde\Xi^{\rm clique}_\ev(\bm \m,\GG)\equiv\m_\ev +\prod_{y\in \obj(\ev)}\Big[1+ \sum_{\ev'\in \fa(y)\atop \ev'\neq \ev}\m_{\ev'}\Big]
\ee
Below we will  refer to the functions $\Xi^{\rm clique}_\ev(\bm \m,\GG)$ and $ \widetilde\Xi^{\rm clique}_\ev(\bm \m,\GG)$
as the {\it clique estimates} of  $\Xi_\ev(\bm \m,\GG)$ and we will make use of  these estimates in Section \ref{appl}.

\\Via the upper bound (\ref{clique}) for $\Xi_\ev(\bm \m,\GG)$,  an improved estimate for latin transversal has been obtained
in \cite{BFPS} and  improved bounds
for  several chromatic indices  have been given in   \cite{NPS} and \cite{BKP}. It is important however to stress that
bounds (\ref{clique}) or (\ref{cliqueb}) are efficient only if the cliques $\{\fa(y)\}_{y\in \obj(\ev)}$ are not too overlapped.
When cliques  $\{\fa(y)\}_{y\in \obj(\ev)}$  are too overlapped these bounds  tends to be too rough. This situation occurs for example
in the case of perfect and separating hash families (see \cite{PS}).

\subsection{Motivation and overview of the present paper}
\\Moser and Tardos brillant insight has represented  a significant advance,
but it is not the end of the story.
Shortly after their work, the question has been raised  as to whether it is possible
to design alternative algorithms, eventually
depending on the specific problem treated, able to go beyond  the
bounds given by the LLL or CELL. These ideas have been
originally developed in \cite{DJKW, GKM, EP}  where   specific algorithms
have been
implemented  for  specific graph coloring problems to obtain
bounds which are  better than those  obtainable  by   LLL or even   CELL.
In particular, in the paper \cite{EP} Esperet and
Parreau  devised an algorithm able to obtain
a new upper bound for the chromatic index  of the acyclic edge coloring  of a graph with maximum degree $\D$ which sensibly
improves on the bound obtained by Ndreca et al.  \cite{NPS}  just an year before via the CELL. Esperet and
Parreau
also suggested  that this algorithm   yields a general method able  to treat most of the
applications in graph coloring problems covered by the LLL.
Indeed, this was confirmed in several successive  papers   \cite{GMP,Pr,OP,PSS,CR, CLY,CR,FS, SV}, where the
Esperet-Parreau scheme has been applied to several
combinarorial problems in graph colorings and beyond, generally improving previous results obtained via the LLL
(sometime the improvement is sensible, sometimes less)
This approach, in some sense a variant of the Moser-Tardos algorithmc  LLL, has been called
{\it entropy compression} method (the expression
was first used by Tao  in reference  to the Moser Tardos algorithm \cite{Ta}).

\\In works mentioned
the entropy compression method is used as  a set of instructions to devise   specific algorithms to face  specific combinatorial problems rather
than a general theorem which
one can use a possible alternative  to  the   LLL.
It is therefore  natural   to ask  whether   there is some general  {\it entropy-compression criterion},
valid  in a (possibly restricted) Moser-Tardos setting
that can be used as an alternative to the  the LLL  and  cluster expansion criteria.
For the specific framework of graph coloring, this question has been  addressed
in a paper  by Gon\c{c}alves, Montassier and Pinlou in \cite{GMP}
where  a general criterion is somehow outlined (see there Theorem 12).
It is also worth  to mention a paper by
Bernshteyn  \cite{Be} in which  a non algorithmic ``Local Cut Lemma", inspired by the Entropy-compression method,  is able to
reproduce some of its  achievements.
A second natural question is to try to understand
to what extent the entropy-compression method is more efficient than the local lemma
and for what reasons.

\\In the present  paper,   basically a survey of entropy compression method, we tried to answer these questions.
We first of all  of clarify the (as general as possible) framework in which the EC-method can be applied.
Namely, the EC-setting coincides with  a restricted  Moser-Tardos variable setting in which the
independent random variables generating
the product  probability space  take  values in a common  finite set of  cardinality $k\in \mathbb{N}$ according to the uniform distribution.
We then show that in such a  setting,  a  general constructive criterion (in Theorem \ref{teo1} below) alternative to the  LLL  and
CELL criteria can indeed be formulated.
This ``Entropy-Compression criterion" (EC criterion), i.e.  inequality (\ref{condentr})  below,
is  in general easy to be implemented
 in applications and, once satisfied, furnishes
a general random algorithm  which  efficiently find  a configuration  avoiding all bad events.

\\We also point out that  in the specific (and, as fas as we know, unique) case of
the acyclic edge coloring, the Entropy-Compression setting  goes somehow beyond this restricted
Moser-Tardos setting in the sense that
bad events not sharing variables are not necessarily mutually independent (and this is an extension
 w.r.t. the Moser-Tardos framework).

\\We finally  show  how to implement the EC criterion  to several well known (and hopefully pedagogical)
problems in combinatorics.
These example will be also used to make a comparison between EC and CE criteria trying to clarify in which
conditions  the EC criterion beats  the  CE and LLL criteria   and to what extent.

\\The rest of the paper is organized as follows. In Section 2 we illustrate  the EC setting  and we  state Theorem \ref{teo1}
containing the general EC criterion.
In Section 3 we give the proof of Theorem \ref{teo1}. In section 4 we discuss some applications.
\def\0{\emptyset}

\section{The entropy-compression  setting}

\\Let $\mathbf{S}$  be finite set and let $k\in \mathbb{N}$. We refer to $\mathbf{S}$ as the set of ``atoms" and to $[k]$ as the set of ``colors"
We suppose that  a total order is fixed  in sets  $\mathbf{S}$ and $[k]$. We let $[k]_0=[k]\cup\{0\}$.
A coloring $\k$  of $\mathbf{S}$ is a function
$\k:\mathbf{S}\to [k]$. A partial coloring $\g$ of $\mathbf{S}$ is a function $\g:\mathbf{S}\to [k]_0$ and
when  $\g(x)=0$ we say that $x$ is uncolored. For any non empty $Y\subseteq \mathbf{S}$ we denote by $[k]^Y$ and $[k]_0^Y$
the sets of colorings and partial colorings
in $Y$ respectively. For any nonempty set $X$, we denote by  ${\mathcal P}(X)$  the set of all nonempty subsets of $X$.

\underline{}\\Given the pair $(\mathbf{S},k)$, let ${\mathcal O}$  be a family of subsets  of $\mathbf{S}$  (called objects).
 Given an object  $A\in \A$,
a {\it  standard flaw of $A$} (the bad event in the Moser-Tardos scheme)   is a subset $\ev\subset [k]^A$. A
a {\it  non-standard flaw of $A$}   is  a function $\mathfrak{f}: [k]_0^{\mathbf{S}\setminus A}\to \PP([k]^A)$ (i.e. for each $\g\in  [k]_0^{\mathbf{S}\setminus A}$, $\ef(\g)$ is
a subset $\ef\subset [k]^A$).
If $\ev$ is a  flaw of $A\in \A$ (either standard or non-standard),  we write $\obj(\ev)=A$ and we call $|A|=|\obj(\ev)|$ {\it the size of $\ev$}
while $|\ev|$ is the cardinality of $\ev$ (i.e. the number of colorings forming $\ev$). As usual, $\bar \ev$ will denote the complement of $\ev$ in $[k]^A$,
i.e. $\bar\ev=[k]^A\setminus\ev$.
A family $\fa$ of flaws is said to be {\it associated to $\A$} if  each $\ev\in \fa$ is a flaw of $A$ for some $A\in \A$.
  When $y\in \obj(\ev)$ we will (improperly)
say that $\ev$ contains $y$ and we denote by $\fa(y)=\{\ev\in \fa: y\in \obj(\ev)\}$ the sub-family of  flaws containing $y$.
A {\it a good coloring w.r.t. to $\fa$} is a coloring $c\in [k]^\mathbf{S}$  avoiding all flaws in $\fa$, that is to say,
$c$ is such that for all $A\in \A$ and all $\ev\in \fa$ such that $\obj(\ev)=A$, we have $c|_A\notin \ev$ if $\ev$ is standard and
$c|_A\notin \ev(c|_{\mathbf{S}\setminus A})$ if $\ev$ in non-standard. In other words, the restriction $c|_A$ of the coloring $c$
to the set $A$ belongs to $\bar\ev$ for each flaw $\ev$ of  $A$.
We denote by $\bigcap_{\ev\in \fa}\bar\ev$ the set of all good colorings w.r.t. $\fa$
in $[k]^\mathbf{S}$ .

\vv
\\{\bf Remark}. When a family $\fa$ associated to $\A$ contains only
standard  flaws (which is what happens in nearly all known applications of the EC method, as far as we know)
the EC-setting  coincides
with a Moser-Tardos variable setting whose  probability space is generated by $|\mathbf{S}|$ i.i.d.  random variables $\{\psi_e\}_{e\in S}$
all taking value in the set $[k]$ according to the uniform distribution. On the other hand
when some of the flaws in the family $\fa$ are non-standard (the only example we know is the case of the acyclic edge coloring),
the EC-setting has no  manifest  correspondence with the
Moser-tardos setting.
We will refer below to first  case   as the {\it independent EC-setting}
and to the second case as the {\it dependent EC-setting}. We stress once again that virtyally all
applications of the entropy compression method available  in the literature fall in the
{\it independent EC-setting}, with the sole exception, far as we know, of the
the acyclic edge coloring of  bounded degree graphs, for  which the use the {\it dependent EC-setting}
is required (see section \ref{appl}).

\begin{defi}\label{reco}
A flaw $\ev$
is said tidy  if, for all $y\in \obj(\ev)$,  there exists   $X\subset \obj(\ev)\setminus \{y\}$
such that any coloring $\k\in \ev$ is uniquely determined by its restriction to $X$
 and no $Y\varsubsetneq X$ has this property. The set $X$ is called  a ``seed"  of   $\ev$.
\end{defi}

\\{\bf Remark}.
 By definition a tidy flaw has  the empty set as its unique seed if and only if is elementary. Note that
if $\ev$ is not tidy, then it can be seen as  the disjont union of tidy (in the worst case  elementary) flaws.
Therefore there is no loss of generality in considering only
families  in which all flaws are tidy.

\vv
\\Given a family $\fa$ of tidy flaws associated to $\A$ and  given a  flaw $\ev\in \fa(y)$,  we  define
\be\label{kya}
\k_y(\ev)= \min\{|X|: ~X\subset \obj(\ev)\setminus \{y\}~ {\rm and}  ~X~{\rm is\;a\; seed\; of\;}\ev\}
\ee
and set
\be\label{power}
\|\ev_{y}\|= |\obj(\ev)|-\k_y(\ev)
\ee
We refer to $\|\ev_{y}\|$ as {\it the  power of the event $\ev$} (w.r.t. $y$).
Note that, for any  $\ev\in \fa(y)$, $\k_y(\ev)\in \{0,1,\dots, |\obj(\ev)|-1\}$ and  $\|\ev_y\|=|\obj(\ev)|$  only if
 $\ev$ is elementary.

\\We set
\def\fa{{\mathfrak{F}}}
$$
\fa_{j}(y)=\{ \ev\in \fa(y)~{\rm and}~\|\ev_{y}\|=j\}, ~~~\fa_j=\bigcup_{y\in \mathbf{S}}\fa_{j}(y)
$$
and
\be\label{dj}
d_j= \max_{y \in \mathbf{S}, \g\in [k]_0^{\mathbf{S}}} \big|\fa_{j}(y)\big|
\ee
I.e.  $d_j$ is an upper bound for the number of events with effective size $j$ associated to objects containing a common atom.
\\We finally define
\be\label{EA}
E_\A=\{j\in \mathbb{N}: \exists y\in \mathbf{S}~  {\rm such ~that }~  \fa_j(y)\neq \0 \}
\ee
and set,
for $x\in (0,+\infty)$,
\be\label{phiea}
\phi_\A(x)=  1+ \sum_{j\in E_\A}d_j x^{j}
\ee

\\{\bf Remark}. Note  that
a tidy event $\ev$
of size 1, i.e. such that $\obj(\ev)=x$, is necessarily elementary because its unique seed is the empty set. Therefore,
give $A\in \A$ such that $|A|=1$ and $\g\in [k]_0^{\mathbf{S}\setminus A}$,
a tidy {\it flaw}  of  $A$ given $\g$ has  always the form
$\ev=\{\mbox{the atom}~ x\in \mathbf{S}~ \mbox{has the color} ~c ~\mbox{and} ~c\in \Phi(x,\g)\}$ where $\Phi(x,\g)$ is some subset $\subset [k]$
of forbidden  colors for the atom $x$ given the partial coloring $\g$ outside $x$.

\vv\vv

%We finally remark that, once introduced the definition of $\k_l$, we have, for free,
%a upper bound for the probability of each event $A$, i.e. if $A\in {{\mathcal O}}^{\k_l}_l$, then
%$P(A)\leq {k^{-(l-\k_l)}}$.

\def\0{\emptyset}
\noindent
%Let ${\mathcal O}$ be a set of events depending on variables $\cal V$.

\subsection{EC-Algorithm and Entropy-Compression Lemma}\zeq

\\We remind that a total order is chosen in the sets  $\mathbf{S}$ and $[k]$.
We also suppose that a total order has been chosen in the sets $\A$ and $\fa$.
We  choose, for each $y\in \mathbf{S}$ and $\ev\in \fa(y)$,  a
unique subset  $G(\ev, y)\subset \obj(\ev)\setminus\{y\}$  such that  $G(\ev,y)$ is a
seed of $\ev$ and $|G(\ev,y)|=\k_y(\ev)$. We also denote
shortly $G^{c}(\ev,y)=\obj(\ev)\setminus G(\ev,y)$. Note that $y\in G^c(\ev,y)$. Given a partial coloring $\g\in [k]_0^{\mathbf{S}}$, given $x\in \mathbf{S}$
  and given a color $s\in [k]$ we denote by $\g]^s_{_x}$
  the partial coloring which coincides with $\g$ in the set $\mathbf{S}\setminus \{x\}$
  and it takes the value
$s$ at  $x$.

\vskip.2cm
\noindent
{\bf EC-Algorithm}.

\\Assume that a set of atoms $\mathbf{S}$ is given together with
a family of objects $\A$ and a family $\fa$ of flaws associated to $\A$.
Let $t$ be an arbitrary natural number (which can be taken as large as we please).
The input of the algorithm is a  vector  $F_t$
in the space  $[k]^t$ (i.e. a vector with $t$ entries each of which taking values in the set $[k]$).

\\The algorithm
performs (at most) $t$ steps and to  each step $i\in [t]$ a partial coloring $\g_i$ is associated as described below.

\vskip.2cm
\\- Step $0$. Set $\g_0=0$ (i.e. in the beginning no atom is colored).
%\\- Step 1- Pick up at random a value  for the first variable (in the order chosen).
%\vskip.2cm
%%

\\- Step $i$ (for $i\ge 1$).
\begin{itemize}
\item [{$i^\circ$}]
If $\g^{-1}_{i-1}(0)\neq\emptyset$,
let $y$ be  the smallest atom (in the total order chosen) left uncolored by the partial coloring $\g_{i-1}$.
Take $i^{th}$ entry of the vector $F_t$ and  let $s\in [k]$ be this entry. Color
$y$ with the  color $s$ and consider the partial coloring $\g_{i-1}]^s_y$ obtained from $\g_{i-1}$ by coloring $y$ with the color $s$.
\begin{enumerate}
\item[$i^\circ_1)$]
If no flaw  occurs given $\g_{i-1}]^s_y$, set
 $\g_i=\g_{i-1}]^s_y$ and go to the step $i+1$.
\item[$i^\circ_2)$]
Conversely, if
some  flaw occur given $\g_{i-1}]^s_y$,
select the smallest, say $\ev$, which by construction
belongs to the set $\fa(y)$.
Set $\g_i=\g_{i-1}]^0_{G^c(\ev,y)}$ and go to the  step $i+1$.
In words $\g_i$ is obtained from $\g_{i-1}$ by  discoloring all atoms  in $\obj(\ev)$
except those belonging to the set $G(\ev, y)$.
\end{enumerate}
\item[{$i^\bullet$}]
If $\g^{-1}_{i-1}(0)=\emptyset$, stops the algorithm discarding all entries $f_i,f_{i+1},\dots,f_t$ of $F_t$.
\end{itemize}
\vskip.2cm
\\
Note that the partial coloring $\g_i$  returned by the algorithm at the end of each step $i$ m necessarily avoids all flaws in $\fa$.
The algorithm performs  at most  $t$ steps but it can stop  earlier, i.e. after having perfomed $m<t$ steps
and $\g_m^{-1}(0)=\emptyset$, (in other words if at the end of step $m$
all variables are colored, and no bad event occurs).
In this case only the first $m$ entries of the vector $F_t$ are used. We say that the
algorithm is successful if it stops after $m<t$ steps,
or it lasts $t$ steps and after the last step $t$ we have $\g_t^{-1}(0)=\emptyset$. Conversely, we say that
the algorithm fails if it performs all $t$ steps and $\g_t^{-1}(0)\neq \0$. Clearly when the EC-algorithm is successful we have found
a coloring $\k\in [k]^\mathbf{S }$ which avoid all flaws of the family $\fa$.
Note that    the algorithm  can be either deterministic if $F_t$  is a given prefixed vector or random
if the entries of $F_t$ are uniformly sampled from the set $[k]$  sequentially and  independently.
\vskip.2cm

 \\We are now in the position to state the main theorem of this note.

 \begin{theo}[{\bf Entropy-compression  Lemma}]\label{teo1} Assume that a  pair $(\mathbf{S},k)$ is given together with
a family of objects $\A$ and a family $\fa$  of flaws associated to $\A$. If
\begin{equation}\label{condentr}
k> \inf_{x>0}{\phi_\A(x)\over x}
\end{equation}

\\then
$$
\bigcap_{\ev\in \fa_\A}\bar \ev\neq\0
$$
Moreover, if inequality (\ref{condentr}) holds stricly,  the (random) EC-algorithm finds a coloring $c\in \bigcap_{\ev\in \fa_\A}\bar \ev$  in an expected number of steps linear in $|\mathbf{S}|$.
\end{theo}

\\We will refer below  inequality (\ref{condentr}) as the {\it entropy-compression criterion}.

\vv
\\{\bf Remark}.  The random EC-algorithm still finds a coloring avoiding  all flaws in the family $\fa_\A$ even if
we only demand
$k$ to be  greater or equal (instead of strictly greater) than the l.h.s. of (\ref{condentr}),
but in this case we have no control on
the expected running time.

\section{Proof of Theorem \ref{teo1}}\zeq
\def\0{\emptyset}\def\EA{{E_\A}}

\\Let ${\cal F}_t$ be the subset of $[k]^t$ formed by all vectors $F_t$ such that $\g^{-1}_t(0)\neq \0$. By definition,
if $F_\t\in {\cal F}_t$, the algorithm performs all the $t$ steps and fails.
Clearly $|{\cal F}_t|\le k^t$ and if we are able to prove that this inequality is strict,
then this means that the set $[k]^t\setminus {\cal F}_t$ is
non empty and for each vector $F_t\in [k]^t\setminus {\cal F}_t$,  the algorithm is succesful.

\vskip.2cm
\\So, let assume from now on that $F_t\in {\cal F}_t$. In this case the algorithm induced by $F_t$ is such that
for any step $i\in [t]$ there is no case $i^\bullet$ and therefore
we may set $i^\circ\equiv i$ for all $i\in [t]$.
We start by observing   that, if at step $i$,  after coloring the atom $y$, we fall in the   $i_2)$ procedure, some occurring  flaw
$\ev$ with effective size  $\|e_y\|=j\in E_\A$ is selected. Since $|\fa_j(y)|\le d_j$,
this event $\ev$  can be  uniquely determined by
a pair  $(j,\ell)$ with  $\ell\in [d_j]$.
%Observe also that at the end of procedure  $i_2)$ the variable $y$ is always uncolored.%

\\The EC-algorithm at step $i$ takes $i^{th}$ entry of the vector $F_t$, say $s\in [k]$, and attributes  to the smallest uncolored
atom at step $i$, say $y$,
the color $s$. The vector $F_t$ determines thus uniquely
all the $t$ steps performed by the EC-algorithm.
We will now show that the converse is also true:
any vector $F_t\in {\cal F}_t$ is  uniquely determined by
a ``record" of the algorithm. This record is specified by  a pair $(R_t, \g_t)$ where $\g_t$ is the partial coloring at step $t$ and
$R_t$ is a string  $r_1\cdots r_i \cdots r_t$  with
$r_i$ being defined as follows.

\begin{itemize}
\item if at step $i$, after coloring the atom, $y$ no bad event occur, set $r_i=0$
 \item if at step $i$,  after coloring the atom $y$, the bad event $\ev\in \fa_j(y)$ is selected,
then it uniquely
determined by the pair  $(j,\ell)$ (where $j\in E_\A$ and $\ell\in [d_j]$)
and  set  $r_i=(j, \ell)$
 .
\end{itemize}

%\\{\it Definition of the  function $C_t$}. Denote by $X_t \subset S$ the set of uncolored atoms
%after the step $t$. If $x\in S\setminus X_t$, then its  color is assigned after $t$ step and
%we call $c_t(x)\in C$ its   assigned color.  Let $c_t: S\to C\cup\{0\}$ be the
%partial coloring that assigns to each atom $x \in S$ either the value $c_t(x)$ if $x \in S\setminus X_t$ or $0$ if
%$x \in X_t$. So the partial coloring $C_t$ tells us which atoms are colored (and at which color) and which are not colored yet
%after $t$ steps of the algorithm.

\vskip.2cm
\\Thus  $F_t\in {\cal F}_t$  determines uniquely  the pair $(R_t, \g_t)$. We call $\cal M$
the map $F_t\mapsto (R_t,\g_t)$.  Let us  prove that $\cal M$ is an injection, which in other words
means that the knowledge of the record $(R_t,\g_t)$
of the algorithm
permits to reconstruct uniquely the vector $F_t$ which had produced that record.

\begin{lem} \label{lema 1}
The set $X_t=\g^{-1}_t(0)$ formed by the uncolored atoms
at the conclusion of  step $t$ is uniquely determined by the string $R_t$.
\end{lem}

\\\textbf{Proof}. We prove the statement by induction on $t$.
The case  $t=1$ is easy.  If $R_1=r_1=0$  at the beginning of step $1$ the first atom $x\in S$  is colored,
no bad event occur  and $x$ stays colored at the end of step 1. So in this case we have
$X_1=\mathbf{S}\setminus\{x\}$.
If $R_1\neq 0$, this means that the coloring the first atom $x\in S$ produces
a bad event which is necessarily  elementary and of effective size equal to one (recall: we are supposing that all events of the family
$\fa$ are tidy). The atom  is thus uncolored at the end of step 1
and so in this case $X_1=\mathbf{S}$.
Suppose now $t \ge 2$. By the induction hypothesis $X_{t-1}$ is uniquely determined by $R_{t-1}$. The smallest atom in the set
$X_{t-1}$,  call it $y$,  is the one that will be colored at step $t$.
If $r_t=0$, this means that the coloring  of the  atom $y$  has produced no bad events and so
$X_t=X_{t-1} \setminus \{y\}$. If $r_t=(j, \ell)$ then  the coloring  of the  atom $y$ has produced
the event $\ev$  uniquely determined by the pair $ (j, \ell)$ among those
of effective size $j$ containing the atom $y$
and we know that all atoms in $\obj(\ev)$ have been uncolored except those in $G(\ev,y)\subset \obj(\ev)$ so that
we have $X_t=X_{t-1} \backslash \{\obj(\ev)\setminus G(\ev,y)\}$. $\Box$

\begin{lem}\label{lema 2} For any $t\in \mathbb{N}$,
the map $\cal M$ assigning to each vector  $F_t\in {\cal F}_t$  the pair $(R_t, \g_t)$ is injective.
\end{lem}

\\\textbf{Proof}: We prove  by induction that the pair $(R_t, \g_t)$ uniquely determines $F_t$ (in other words we prove that
$|{\cal M}^{-1}(R_t, \g_t)|=1$). Let us start by considering the case  $t=1$. Suppose first $R_1=r_1=0$.  Then the
first atom $x\in \mathbf{S}$ is colored with the colour say $s\in [k]$, no bad event occurs and
 the partial coloring $\g_1$ is such that $\g_1(x)=s$  while
 $\g_1(y)=0$ for all $y\in \mathbf{S}\setminus\{x\}$. Then  the (one-dimensional) vector $F_1=(f_1)$ has entry $f_1=s$.
 On the other hand,
if $r_1=(1,\ell)$ then necessarily a tidy elementary event $\ev$ such that $|\obj(\ev)|=\|\ev_x\|=1$ is selected  with $\ev$
being the event that
$x$ receive  the color, say $s$ which is the   $\ell^{\rm th}$ color of the list $\Phi(x,\g_0)$
and thus the (one-dimensional) vector $F_1=(f_1)$
has entry $f_1=s$.
Suppose now the claim true for  $t-1$.
Namely, given the pair $(R_{t-1},\g_{t-1})$,  we know the vector $F_{t-1}$.
Therefore we have to show that, given the pair $(R_{t},\g_{t})$ we must be able to find the entry $f_t$
of the vector $F_t$ and the function $\g_{t-1}$
(so that, by induction, $(R_{t-1},\g_{t-1})$ gives us all the remaining entries of the vector $F_t$).
By Lemma \ref{lema 1} the knowledge of $R_{t-1}$ determine uniquely  the set $X_{t-1}$  formed by the uncolored atoms
at the conclusion of  step $t-1$  and hence
we know the smallest atom uncolored after $t-1$ steps; call it $y$. We have to consider two cases: a) $(R_t,\g_t)$ is such that $r_t=0$;
b) $(R_t,\g_t)$ is such that $r_t=(j,\ell)$. The case a) is easy. Indeed if $r_t=0$ then the variable $y$ will be colored  at step $t$
and so if $\g_t(y)=s$, then  $f_t=s$ while $\g_{t-1}$ is such that
$\g_{t-1}(x)=\g_t(x)$ for all  $x\in \mathbf{S}\setminus\{y\}$ and $\g_{t-1}(y)=0$.
Let us now consider the case b) i.e. $r_t=(j,\ell)$ where $j\in E_\A$ and $\ell\in [d_j]$. The pair
$(j,\ell)$ uniquely determines the  event $\ev\in \fa_j(y)$ that is selected among those occurring
at the beginning of step $t$ after the coloring of the variable $y$.
Concerning  this event $\ev$,
we know the  subset
$G(\ev,y)\subset \obj(\ev)$  constituted by the atoms  that will continue to stay colored after the conclusion of step $t$
(recall: $G(\ev,y)$ is empty if $\ev$ is elementary).  Now we have $\g_t$
so we know the colors of the atoms (eventually) belonging to $G(\ev,y)$. Since $G(\ev,y)$ is a seed of $\ev$,
the knowledge of the colors of the atoms of   $G(\ev,y)$ (or the fact that $\ev$ is elementary if $G(\ev,y)=\0$)
allows us to deduce which were the colors  of all the other
atoms in $G^c(\ev,y)=\obj(\ev)\setminus G(\ev,y)$. Denote by   $c(x)$  the color of the atom $x\in G^c(\ev,y)$  uniquely determined
by the coloring of $G(\ev,y)$.  Then  $\g_{t-1}(x)=c(x)$ for   $x\in  G^c(\ev,y)\setminus \{y\}$
and $\g_{t-1}(y)=\g_t(y)=0$. Finally, since by definition $y\in G^c(\ev,y)$, if $c(y)$ is, say, equal to $s\in [k]$
then  $f_t=s$. $\Box$
\vskip.2cm

\\Let us denote, for $r\in [|\mathbf{S}|]$,  ${\cal F}^r_t$ the set formed by all vectors in ${\cal F}_t$ such that
$|\g^{-1}_t(0)|=r$.
Clearly ${\cal F}_t$ is the disjoint union of the family $\{{\cal F}^r _t\}_{r\in [|\mathbf{S}|]}$ and therefore we have that
\be\label{sos}
|{\cal F}_t|=\sum_{r=1}^{|\mathbf{S}|}|{\cal F}^r _t|
\ee

\\Let $({\cal R}_t,{\cal G}_t)$ (resp.  $({\cal R}^r_t,{\cal G}_t^r)$) set of all records
 produced with vectors in ${\cal F}_t$ (resp. ${\cal F}^r_t$),
in other words $({\cal R}_t,{\cal G}_t)={\cal M}({\cal F}_t)$ (resp.$({\cal R}^r_t,{\cal G}^r_t)={\cal M}({\cal F}^r_t)$).
%We have, as a consequence of Lemma  \ref{lema 2}, the following proposition.
\begin{prop}
\be\label {pro}
|{\cal F}_t|\le (|[k]|+1)^{|\mathbf{S}|} \sum_{r=1}^N|{\cal R}^r_t|
\ee
\end{prop}
{\bf Proof}.
From
Lemma \ref{lema 2} it follows that
\be\label{sas}
|{\cal F}^r _t|\le |{\cal R}^r _t| |{\cal G}^r _t|\le |{\cal R}^r _t| |{\cal G} _t|\le |{\cal R}^r _t| (k+1)^{|S|}
\ee
since  ${\cal G} _t$
is a subset of the set of all partial colorings from $\mathbf{S}$ to $[k]\cup\{0\}$ which has cardinality $(k+1)^{|\mathbf{S}|}$. Putting
(\ref{sas}) into (\ref{sos}) inequality (\ref{pro}) follows. $\Box$

\vskip.2cm
\subsection{Upper bound for $|{\cal R}^r_t|$}
A string $w_1w_2\dots w_n$ with $w_i\in \{0,1\}$ is  usually called
a word on the alphabet $\{0,1\}$. Give a word $w$ on the alphabet $\{0,1\}$, the
mirror image $\widetilde{w}$ of $w$ is the word obtained by reading $w$ from right to left, i.e
if $w=w_1w_2\dots w_n$ then $\widetilde{w}=w_nw_{n-1}\dots w_1$

\\An initial segment of the word $w_1\dots w_n$
is the a sub-string of  $w_1\dots w_n$ of  the form  $w_1\dots w_i$ with $1\le  i\le n$.
A {\it partial Dyck word} is word on the alphabet $\{0,1\}$ such that in any initial segment of the word the number
of $0$'s is greater or equal than the number of $1$'s. A {\it Dyck word} on the alphabet $\{0,1\}$ is
a partial Dyck word with equal number of $0$'s  and $1$'s (hence a Dyck word has always an even number of letters).
A partial Dyck word can be viewed as a path in $\mathbb{Z}^2$ starting at the origin made by steps either $(1,1)$ or $(1,-1)$
in such a way that the path stays in the first quadrant (i.e. the path never goes below the $x$ axis).
A Dyck word (of size $2n$) is then a Dyck path which starts at the origin and ends at the   point $(2n,0)$ of the $x$ axis.
A descent (ascent) in a partial Dyck word  is a sequence of $1$'s ($0$'s).
For ${E}\subset \mathbb{N}$ let us denote by ${\cal D}_{n,E}$ the set of  Dyck words
with $n$ $0$'s, $n$ $1$'s and with descents having cardinality in $E$. Moreover, for $r<n$ natural number, let us denote
by ${\cal D}_{n,r,E}$  the set of partial Dyck words
with $n$ $0$'s, $n-r$ $1$'s and with descents having cardinality in $E$.

\begin{lem}\label{esppa1}
Let $m_E=\min\{E\setminus\{1\}\}$.  Then the following inequality holds.
\be\label{dick}
|{\cal D}_{n,r,E}|\le |{\cal D}_{t+r(m_E-1),E}|
\ee
\end{lem}
{\bf Proof}.
We can construct an injective function  $f:  {\cal D}_{n,r,E}\to {\cal D}_{n+r(m_E-1),E}$
by associating  to  the partial Dick word $w\in{\cal D}_{t,r,E}$  the Dick word
$f(w)\in {\cal D}_{t+r(m_E-1),E}$ formed by appending at the end of $\a$ the word $(0^{m_E-1}1^{m_E})^r$. $\Box$
\vv
\\Let $w\in {\cal D}_{n,E}$ and let $j\in E$, we define $\d_j(w)$  as the number of descents  in $w$ having cardinality $j$. We further define
$\a_j(w)$  as the number of ascents   in $w$ having cardinality $j$.
Let $\T_{n+1, E}$ be the set of plane rooted trees with $n+1$ vertices such that the number of children of each internal vertx takes value in $E$.
If $\t\in \T_{n+1, E}$ we define $\d_j(\t)$ as the number of  vertices in $w$ having $j$ children.

\begin{lem}\label{esppa2}
There is  a one-to-one map $h: {\cal{T}}_{n+1, E}\to {\cal D}_{n,E}: \t\mapsto w$ such that for any  $j\in E$
\be\label{lemma7}
\d_j(\t)= \d_j(w)
\ee
\end{lem}
{\bf Proof}. We construct the map $h$ as follows. For a given $\t\in {\cal{T}}_{n+1, E}$,
consider the string $\r_t= v_1\dots v_n$ formed by the first $n$ vertices of $\t$
ordered according to a Deep-first search (these are all vertices of $\t$ but a leaf). Now substitute each  $v_i$ in the string $\r_t$
with  the word
$$
\overbrace{0 0\dots 0}^{s_i~{\rm times}}1
$$
where $s_i$ is the number of children of the vertex $v_i$. Let $u$ the word so obtained. By observing that
$\sum_{i=1}^n s_i=n$ we conclude that $u$ is a Dyck word of size $2n$ whose ascents have cardinality in $E$ and
$\a_j(u)=\d_j(\t)$  for all $j\in E$. Moreover the map $\t\mapsto u$  is clearly a bijection.
Now consider the mirror image $\widetilde{u}$ of $u$ and interchange the $0$'s and the $1$'s and let $w$ the word so obtained.
The composition of maps  $\t\to u\mapsto \widetilde{u}\mapsto w$
is clearly a bijection from ${\cal{T}}_{n+1, E}$ to ${\cal D}_{n, E}$ and the word $w$ is a Dyck word of size $2n$ whose descents  have cardinality in $E$ and,
for all  $j\in E$, it holds that $\d_j(w)=\d_j(\t) $. $\Box$

\vv
\\We now construct a (non injective) map $ \mathfrak{M}: {\cal R}^r_t \to {\cal D}_{t,r,\EA}$.
\begin{defi}
\\Let $R_t=r_1\dots r_t\in {\cal R}^r_t$. We define the map $ \mathfrak{M}$ which  associates to $R_t=r_1\dots r_t$ the word
$R^*_t=r^*_1\dots r^*_t$ on the alphabet $\{0,1\}$
 as follows. If $r_i=0$ then put $r^*_i=0$. If $r_i=(j,\ell)$ then put
$$
r^*_i=0\overbrace{1 1\dots1}^{j~{\rm times}}
$$
\end{defi}
\begin{lem}\label{surmap}
For any $R_t\in {\cal R}^r_t$, the word $R^*_t$ defined above
is a partial Dyck word with $t$ $0$'s and $t-r$ $1$'s and descents
in the set $E_\A$ defined in (\ref{EA}).
In other words if  $R_t\in {\cal R}^r_t$, then $R^*_t\in {\cal D}_{t,r,E_\A}$.
\end{lem}
{\bf Proof}. Reading $R^*_t$ from left to the right, every $0$ correspond to a variable that has been colored  and every $1$ correspond
to a variable that has been uncolored if some bad event occurs. So the number of $0$'s is $t$. Let  $x$ be the number of $1$'s. This number
represents by construction the total number of variables that have been uncolored during the process of the algorithm. Thus
$t-x=r$, i.e. $x=t-r$. The string is by construction a
partial Dyck word since we cannot uncolor more variables than the number of colored variable. Finally again by construction we have that the
descents are of size $j\in \EA$.

\begin{lem}\label{surmap2}
The pre-image of a string $R^*_t$ by the map $ \mathfrak{M}$  has cardinality less or equal than
$$
\prod_{j\in \EA}(d_j)^{\d_j(R^*_t)}
$$
where recall that $\d_j(R^*_t)$ is the number of descents in  $R^*_t$ of size $j$ and $d_j$
is defined in (\ref{dj}).
\end{lem}
\\{\bf Proof.}
For each descent of size $j$ in $R^*_t$, we do not distinguish which
flaw of size $j$ that generated this descent. Since there are at most $d_j$ flaws of size $j$ containing a fixed atom,
this accounts for the $\prod_{j\in \EA}(d_j)^{\d_j(R^*_t)}$ factor.
$\Box$

\vskip.2cm
%\\{\bf Remark}. Note that
%\be\label{efi}
%\phi_E(x)= 1+ \sum_{l \in E}\sum_{\k\in K_l} x^{l-\k}= 1+\sum_{j\in F} m_j x^j
%\ee

\\Let us now define
\be\label{fij}
\varphi_j=\cases{1 & if $j=0$\cr\cr
 d_j &if $j\in E_\A$\cr\cr
0 & otherwise
}
\ee
Note that, for $x>0$
\be\label{varp}
\sum_{j\ge 0} \varphi_jx^j= \phi_\A(x)
\ee
We now can prove the following proposition.
\begin{prop}
Let ${\cal T}_{n}$ denotes the set of all plane rooted trees with $n$ vertices, then we have
\be\label{ineq}
|{\cal R}^r _t|\le \sum_{\t\in {\cal T}_{t+r(m_{E_\A}-1)+1}}
\prod_{j\ge 0}\varphi_j^{\d_j(\t)}
\ee
where recall that $\d_j(\t)$ denotes the number of vertices with $j$ children in  $\t$.

\end{prop}
\\{\bf Proof}. From Lemmas \ref{surmap} and \ref{surmap2} it follows that
$$
 |{\cal R}^r _t|\le \sum_{w\in {\cal D}_{t,r,\EA}}\prod_{j\in E_\A}(d_j)^{\d_j(w)}
$$
Lemma \ref{esppa1} then implies
$$
\sum_{w\in {\cal D}_{t,r,\EA}}\prod_{j\in E_\A}(d_j)^{\d_j(w)}\le
\sum_{w\in {\cal D}_{t+r(m_{E_\A}-1),E_\A}}\prod_{j\in E_\A}(d_j)^{\d_j(w)}
$$
Lemma \ref{esppa2} yields
$$
\sum_{w\in {\cal D}_{t+r(m_{E_\A}-1),E_\A}}\prod_{j\in E_\A}(d_j)^{\d_j(w)}=
\sum_{\t\in {\cal T}_{t+r(m_{E_\A}-1)+1,E_\A}}\prod_{j\in E_\A}(d_j)^{\d_j(\t)}
$$
and finally using  definition (\ref{fij}) we get
$$
\sum_{\t\in {\cal T}_{t+r(m_{E_\A}-1)+1,E_\A}}\prod_{j\in E_\A}(d_j)^{\d_j(\t)}=
 \sum_{\t\in {\cal T}_{t+r(m_{E_\A}-1)+1}}\prod_{j\ge 0}\varphi_j^{\d_j(\t)}
$$
$\Box$
\vv
\\We can now  use the following well known theorem in analytic combinatorics  (see e.g., Theorem 5 in \cite{Dm}).
\begin{theo}\label{dmro}
Let $\{\varphi_j\}_{j\ge0}$ a sequence of nonnegative numbers and  let  $\varphi(x)= \sum_{j\ge 0} \varphi_j x^j$; let $R$ be the convergence radius of $\varphi(x)$ and let
$\x$ be the first root in $[0,R)$ of the equation $x\varphi'(x)-\varphi(x)=0$.
%Set $d={\rm gcd}\{j>0:~ \varphi_j>0\}$.
Then
\be\label{dro}
\sum_{\t\in {\cal T}_{n}}\prod_{j\ge 0}\varphi_j^{\d_j(\t)}
%d\left[\varphi(\x)\over 2\pi \varphi''(\x)\right]{[\varphi'(\x)]^n\over n^{3/2}}(1+ {\cal O}(n^{-1})
\le C_{\varphi}
 {[\varphi'(\x)]^n\over n^{3/2}}
\ee
with $C_{\varphi}$ constant.
\end{theo}
Note that
$$
\varphi'(\x)=\min_{x>0}{\varphi(x)\over x}
$$
Inserting now inequality (\ref{dro}) into (\ref{ineq}) and recalling  (\ref{varp}) we get
\be\label{ineq2}
|{\cal R}^r _t|\le
 C_{\varphi} {[\phi'_\A(\x)]^{t+r(m_{E_\A}-1)+1}\over ({t+r(m_{E_\A}-1)+1})^{3/2}}
\ee

\subsection{Conclusion of the proof of Theorem \ref{teo1}}\label{conclu}
We are now in the position to end the proof of Theorem \ref{teo1}. Indeed,
inserting (\ref{ineq2}) into
(\ref{pro}) we get
$$
|{\cal F} _t|~\le ~  C_{\varphi} (k+1)^{|\mathbf{S}|} \sum_{r=1}^{|\mathbf{S}|}
 {[\phi'_\A(\x)]^{t+r(m_{E_\A}-1)+1}\over ({t+r(m_{E_\A}-1)+1})^{3/2}}~\le
$$
$$
\le  C_{\varphi,\mathbf{S}, k} \left[\phi_\A'(\x)\right]^{t}
{1\over {t}^{3/2}}
$$
with
$$
C_{\varphi,\mathbf{S}, k}= C_{\varphi} (k+1)^{|\mathbf{S}|} [\phi_\A'(\x)]^{{|\mathbf{S}|}(m_{E_\A}-1)+1}
$$
Hence
if
\be\label{opa}
\left[ \phi_\A'(\x)\right]=\min_{x>0}{\phi_\A(x)\over x}\le k
\ee
we have that for $t$  larger than $C_{\varphi,\mathbf{S}, k}^{2/3}$
$$
|{\cal F} _t|< k^t
$$
i.e. the algorithm stops. The condition (\ref{opa}) which ensures
the algorithm to stop is exactly (\ref{condentr}) with greater or equal replacing strictly greater.
Now note that, if each entry of the vector ${\cal F} _t$ is selected uniformly and independently at random, then  $|{\cal F} _t|/k^t$ is the probability that after $t$ steps the (random) EC-algorithm is still running.
Let
us  suppose that inequality (\ref{opa}) holds strictly. I.e., let us suppose that there exists $\e>0$ such
that $({ \phi_\A'(\x)/k})\le (1-\e)$.
Then
the probability $P(t)$
that the algorithm is still running  after $t$ steps can be bounded above as follows
$$
P(t)\le C_{\varphi} (k+1)^{|\mathbf{S}|} [\phi_E'(\x)]^{|\mathbf{S}|(s-1)+1}{(1-\e)^t}
$$
Let $t_0$ be the solution of the equation
$$
C_{\varphi} (k+1)^{|\mathbf{S}|}N [\phi_E'(\x)]^{|\mathbf{S}|(s-1)+1}{(1-\e)^t}=1
$$
Remark that $t_0$ is linear in $|\mathbf{S}|$ and observe that
$P(t_0+t)\le (1-\e)^t=e^{-t|\ln(1-\e)|}$.
Therefore, if  inequality (\ref{opa}) holds strictly, the expected running  time $T$ of the EC-algorithm is bounded by
$$
T\le t_0+{1\over |\ln(1-\e)|}
$$

%\vskip4.0cm
%Note that $|{\cal F} _t|/k^t$ is the probability that after $t$ steps the algorithm is still running.
%Let
%us now suppose that inequality (\ref{opa}) holds strictly. I.e., let us suppose that there exists $\e>0$
%such that $({ \phi_E'(\t)/k})\cdot\max_{l,\k} (d_l^{\k})^{1/(l-\k)}\le (1-\e)$.
%Then
%the probability
%that the algorithm stops after $t$ steps is at least ${1\over 2}$  as soon as  $|{\cal F} _t|/k^t<1/2$. This yields
%$$
%C_{\varphi} (k+1)^N [\phi_E'(\t)]^{N(s-1)+1}{(1-\e)^t\over t^{3/2}}<
%{1\over 2}
%$$
%which is surely  satisfied if
%$$
%t>  {N\over \e}\ln\left[2\tilde C_{\varphi} (k+1)[\phi_E'(\t)]^{s}\right]
%$$
%where we have set $\tilde C_\varphi=\max\{1, C_{\varphi} ^{1\over N}\}$.
\\So if the condition  (\ref{opa}) is replaced by the  inequality
\be\label{opa2}
 \phi_\A'(\x)\le (1-\e)k
\ee
for some $\e>0$,
then the expected running time of the algorithm is linear in $|\mathbf{S}|$.

\\This concludes
the proof of Theorem \ref{teo1}.

\section{Comparison and examples}\label{appl}\zeq

\subsection{Entropy-Compression   versus cluster expansion}\label{compar}
It is instructive and illuminating, we think,  to  start by
comparing the entropy-compression criterion with the cluster expansion criterion on a general basis.
\\We will assume here that we are in the  {\it independent} EC-setting.
For a  comparison between EC and CELL in the {\it dependent}  EC-setting (i.e. the case of the acyclic edge coloring) we refer  the reader to Section \ref{subacy}.

\\We have therefore
the pair $(\mathbf{S},k)$  together with
a family of objects $\A$ and a  family $\fa$ of standard flaws associated to $\A$.
As previously mentioned, this  corresponds to a Moser-Tardos setting
whose  probability space is generated by $|\mathbf{S}|$ i.i.d.  random variables $\{\psi_x\}_{x\in \mathbf{S}}$
all taking value in $[k]$ according to the uniform distribution.
In this framework
it is natural (and easier) to use  the clique estimate   (\ref{clique}) in order to evaluate  $\Xi_\ev(\bm \m,\GG)$. Therefore the CE criterion
(\ref{cond1n}) is satisfied if
\be
\label{condcli}
\prob(\ev)\;\leq\;  {\mu_\ev \over\Xi^{\rm clique}_\ev(\bm \m,\GG)}\;=\; {\mu_\ev \over  \prod_{y\in \obj(\ev)}\Big[1+ \sum_{\ev'\in \fa(y)}\m_{\ev'}\Big]}
\ee
We  will refer  below to condition (\ref{condcli}) as the {\it strong variable setting CE criterion}.
This condition (\ref{condcli}), of course, is strongest than the CE criterion (\ref{cond1n}) in the sense that to be fulfilled
requires smaller probabililies for the bad events.
We stress once again that   the {strong variable setting CE criterion} (\ref{condcli})
 has been the tool used in recent times to improve several bounds previously obtained via the classical version of LLL given by Theorem \ref{lov}.

\\Now,  in the  elementary probability space generated by i.i.d. uniformly distributed random variables  $\{\psi_x\}_{x\in \mathbf{S}}$
taking values in $[k]$,  given  $y\in \mathbf{S}$, $j\in E_\A$, and  $\ev\in \fa_j(y)$ and recalling definitions (\ref{kya}) and (\ref{power}),
it is immediate to see that
\be\label{probb}
\prob(\ev)= {k^{\k_y(\ev)}\over k^{|\obj(\ev)|}}= {1\over k^{|\obj(\ev)|-\k_y(\ev)}}  =  {1\over k^{j}}
\ee
Moreover, for any $\ev\in\fa$
\begin{equation}\label{des}
\Xi^{\rm clique}_\ev(\bm \m,\GG)=\prod_{y\in \obj(\ev)}(1+ \sum_{\ev'\in \fa(y)}\m_{\ev'})\le
\prod_{y\in \obj(\ev)}(1+ \sum_{j\in E_\A}\sum_{\ev'\in \fa_j(y)}\m_{\ev'})
\end{equation}
Given   (\ref{probb}), it is  natural to set $\m_\ev= x^j$   for all $\ev\in \fa_j$ , where $x>0$, so that,
recalling definition (\ref{dj}),  we get
\be\label{cli2}
\Xi^{\rm clique}_\ev(\bm \m)\le
\prod_{y\in \obj(\ev)}(1+ \sum_{j\in E_\A}d_jx^j)= \Big[\phi_\A(x)\Big]^{|\obj(\ev)|}\le  [\phi_\A(x)]^{l_j}
\ee
where
$$
l_j=\max\{|\obj(\ev)|: \ev\in\fa_j\}
$$
Note that $\phi_\A(x)$ is the very same function appearing in the entropy-compression criterion
(\ref{condentr}). Using (\ref{cli2}), condition (\ref{condcli}) is rewritten as
$$
{1\over k^{j}}\le\max_{x>0}{x^{j}\over  [\phi_{\A}(x)]^{l_j}} ~~~~~~~\mbox{ for  all~~~ $j\in E_\A$}
$$
In conclusion, in this restricted Moser-Tardos variable setting in which the
random variables  $\{\psi_x\}_{x\in \mathbf{S}}$ are i.i.d. uniformly distributed all taking values in the same set $[k]$,
the strong variable setting CE criterion (\ref{condcli}) can be rewritten as
\be\label{opa3}
k\ge \min_{x>0}{\left[\phi_\A(x)\right]^{l_j\over j}\over x}~~~~~~~\mbox{ for  all~~~ $j\in E_\A$}
\ee
It is now crystal clear to compare (\ref{opa3}) with the  entropy-compression criterion  (\ref{condentr}).

\\Since for all $j\in E_\A$  we have that $l_j\ge j$ and  $\phi_\A(x)>1$,  it holds
$$
{\left[\phi_\A(x)\right]^{l_j\over j}\over x}\ge  {\phi_\A(x)\over x}~~~~~~~\mbox{ for  all~~~ $j\in E_\A$ and for all  $x>0$}
$$
%and hence, for all  $j\in E'_\A$ , we have
%\be\label{okkk}
%\inf_{x>0}{\left[\phi_\A(x)\right]^{l_j\over j}\over x}~\ge \inf_{x>0}{\left[\phi_\A(x)\right]\over x}
%\ee
and the equality only holds if $l_j=j$ for all $j\in E_\A$ which happens only when all (tidy) events of tha family $\fa$ are elementary.

\\The conclusion is that  in the independet EC-setting $(\mathbf{S},k)$, the entropy compression criterion (\ref{condentr}) is always better than  or equal to condition (\ref{opa3})
and these two conditions coincide  only all  events of tha family $\fa$ are elementary. However, it is worth to remind
that condition (\ref{opa3}) has been derived from (\ref{condcli}) where the clique estimate
$\Xi^{\rm clique}_\ev(\bm \m,\GG)$ has been used. As said above, condition (\ref{condcli})
yields bounds poorer  than those  obtainable by the  general CE criterion (\ref{cond1n}) as soon as the cliques composing $\G_\GG(\ev)$ are
too overlapped.
%We will give in section (\ref{PropB}) and example about this circumstance.

\subsection{Some graph coloring examples: when and to what extent EC criterion beats CE criterion}
In this subsection $G=(V,E)$ is  a graph with maximum degree $\D$  and  $k\in \mathbb{N}$.
\subsubsection{Non repetitive vertex coloring}

Let us start with a simple example fitting in the {\it independent} EC-setting.
A coloring $c: V\to [k]$ of the vertices of $G$ is called {\it nonrepetitive} if, for any $n \ge 1$, no path $p=\{v_1,v_2,\dots, v_{2n}\}\subset V$
is repetitive (i.e such that $c(v_i)=c(v_{i+n})$ for all $i=1,2,\dots n$). Observe that a nonrepetitive vertex
coloring is proper.

\noindent
The  EC-setting $(\mathbf{S},k)$ for the present case is as follows. The set of atoms is $\mathbf{S}=V$  and  the set of colors is $[k]$.
We now determine
the families of objects $\A$ and of  flaws $\fa$ associated to $\A$.
Let $P_{n}$ be the set of all paths with $2n$ vertices in $G$ and let $P=\cup_{n \geq 1}{P_n}$.
We set  $\A= P$. Let now for $p\in P$, $\ev_p$ be the event  that $p$ is monochromatic and set
$
\fa_P=\{\ev_p\}_{p \in P}
$
Clearly any event of the family $\fa_P$ is standard and therefore here, as in all the examples to follow (except that
of Subsection (\ref{subacy})), we are in the {\it independent} EC setting. Moreover
observe that any event $\ev\in \fa_P$ is  tidy and, for any $v\in \obj(\ev)$,  we have that $\|\ev_v\|={1\over 2}|\obj(\ev)|$.
Therefore in the present case  we have

\noindent
$$
E_\A= \{j\}_{j\ge 1},
$$
Furthermore, for $j\in  E_\A$, the number $d_j$ defined in (\ref{dj}) represents here
the maximal number of paths of size $2j$ in $G$ containing a fixed vertex, which can be bounded as
\be\label{djei}
d_{j} \leq j\D^{2j-1}
\ee
so that
$$
\phi_\A(x)\le 1+ \sum_{j=1}^\infty j\D^{2j-1}x^j,  ~~~~ \mbox{with $0<\D^2x<1$}
$$
Hence the EC criterion  (\ref{condentr}) in the present case is fulfilled if
\be\label{repve}
k\;\ge \;\min_{0<\D^{2}x<1}{1+{1\over \D}\sum_{j\ge 1}j(\D^{2}x)^j\over x}
\ee
Inequality (\ref{repve}) can be rewritten as
$$
{k}\ge \Big[\min_{0<x<1}g(x)\Big]\D^2
$$
where
$$
g(x)={1+ {1\over \D}{x\over (1-x)^2}\over x} = {1\over x}+{1\over \D}{1\over (1-x)^2}
$$
%The minumum of $g(x)$ in the interval $x\in (0,1)$ occurs in $x_0$ where $x_0$ is the solution of the equation
%$$
%\D(1-x)^3-2x^2=0
%$$
 For $\D\ge 3$
we can take $x_1= 1-2^{1\over 3}/\D^{1\over 3}$ (which  is such that   $x_1-\min_{0<x<1}g(x)= o(\D^{-1/3})$)
and  we can bound
$$
g(x_0)\le g(x_1)=
 1+ {3\over 2^{1\over 3} \D^{1\over 3}}+ {2^{2\over 3}\over \D^{2\over 3}}{1\over 1- {2^{1\over 3}\over \D^{1\over 3}}}, ~~
 $$
Hence  (\ref{repve}) is satisfied if
\be\label{acyc}
k\ge \D^2\Big( 1+ {3\over 2^{1\over 3} }\D^{-{1\over 3}}+ o({\D^{-{1\over 3}}}) \Big)
\ee
which is
the  bound   given  in \cite{GMP}.

\noindent
It is worth to notice that in this first example the strong variable setting CE criterion in the form (\ref{opa3})
obtains a bound
nearly as good as (\ref{acyc}). Indeed, recalling that  in the present case $\|\ev_v\|={1\over 2}|\obj(\ev)|$
for any  $\ev\in \fa$  and any $v\in \obj(\ev)$,   we have that ${l_j\over j}=2$
for all $j\in E_\A$.  Therefore
(\ref{opa3}) is written as
\be\label{clus}
{k}\ge \Big[\min_{0<x<1}\widetilde g(x)\Big]\D^2
\ee
with
$$
\widetilde g(x)={\big(1+ {1\over \D}{x\over (1-x)^2}\big)^2\over x}
$$
and, by numerical computation (\ref{clus}) is satisfied if
$$
k\ge \D^2\Big( 1+ {5\over 2} \D^{-{1\over 3}}+o(\D^{-{1\over 3}})\Big)
$$
which is  worst than the   bound  (\ref{acyc}) obtained via EC criterion only at the order $O(\D^{-{1\over 3})}$.

\subsubsection{Star vertex coloring}
A proper coloring $c: V\to [k]$ of the vertices of $G$ is called
a {\it star coloring} if no path $p=\{v_1,v_2,v_3,v_4\}$
is bichromatic. The mimumum numbers of colors required
for a graph $G$ to have a vertex star coloring is called the star chromatic index of $G$ and  it is denoted by $\pi(G)$.

\\Let $P_2$   be the set of all paths with $2$ vertices and let $P_4$ be the set of all paths with $4$ vertices.
The EC-setting in the the present case  is established by taking $\mathbf{S}=V$, $\A=P_2\cup P_4$. The family $\fa_{P_2\cup P_4}$ of flaws
(a subfamily of the family $\fa_P$ of the previous example) is  given by
$\fa= \fa_{P_2}\cup\fa_{P_4}$
where
$\fa_{P_2}= \{\ev_p\}_{p\in P_2}$ and $\fa_4 = \{\ev_p\}_{p\in P_4}$.
As in the previous example all flaws of in $\fa_{P_2}\cup\fa_{P_4}$ are standard,  non elementary and  tidy. This is once again an {\it independent} EC-setting.
Clearly, a coloring $c: V\to [k]$  avoiding all   flaws in $\fa$ is a star coloring.
 Moreover $\ev\in \fa_{i}$ ($i=2,4$) is
such that  for any $v\in \obj(\ev)$,  it holds  $\|\ev_v\|={i\over 2}$.
This implies that
$$
E_\A= \{1,2\}
$$
and, recalling (\ref{djei})
$$
\phi_\A(x)\le 1+\D x+2\D^{3}x^2,  ~~~~ \mbox{with $x>0$}
$$
Hence the entropy-compression criterion  (\ref{condentr}) in the actual case is fulfilled if
\be\label{repve2}
k\ge \min_{x>0}{1+\D x+2(\D^{3}x^2)\over x}
\ee
with $x>0$.
Setting $\D x=y$, inequality (\ref{repve2}) can be rewritten as
\be\label{estar}
{k}\ge \left(\Big[\min_{y>0}g(y)\Big]\right)\D
\ee
where
$$
g(y)={1+y+ 2\D y^2\over y}
$$
The minumum of $g(y)$  for $y>0$  occurs at $y_0=(\sqrt{2}\D^{1\over 2})^{-1}$ and $g(y_0)=\sqrt{2}\D^{1\over 2}(2+(\sqrt{2}\D^{1\over 2})^{-1})$
and thus
\be\label{astar}
k\ge 2\sqrt{2}\,\D^{3\over 2}+\D
\ee
which is the  bound given in \cite{EP}.

\\Finally, let us compare the bound (\ref{estar}) with that produced by the strong variable setting
CE criterion (\ref{opa3}). As in the previous example, ${l_j\over j}=2$ for $j=1,2$, so that  (\ref{opa3}) is
\be\label{clstar}
{k}\ge \left(\Big[\min_{y>0}{(1+y+ 2\D y^2)^2\over y}\Big]\right)\D
\ee
which yields
$$
k\ge {8\over\sqrt{27}}\,2\sqrt{2}\,{\D^{3\over 2}} +{8\D\over 3}+ o(\D)
$$
being worst  than  (\ref{astar}) by a factor 1.54 at the lowest order in $\D$.

\subsubsection{\bf Vertex coloring graphs frugally}
Given a natural number $\b$,
a proper coloring of the vertices of $G$ is  {\it $\b$-frugal} if any vertex has at most
$\b$  members of any color class in its neighborhood.
The minimum number of colors required such that a graph $G$ has
at least one $\b$-frugal proper vertex coloring is called the {\it $\b$-frugal chromatic number} of $G$ and
will be denoted by $\chi_\b(G)$.

\\Given $v\in V$, we denote by
$H_{\b}^v$ the set of all subsets of the neighborhood $\G_G^*(v)$ of $v$ with exactly $\b+1$ vertices. In other words
$$
H_\b^v=\{Y\subset \G^*(v): ~|Y|=\b+1\}
$$
Moreover, given $v\in V$ and $h_\b^v\in H_\b^v$, we let $\ev_{h_\b^v}$ to be the event that
$h_\b^v$ is monochromatic).
Let ${\cal H}_\b=\{H_\b^v\}_{v\in V}$ and let $\fa_{{\cal H}_\b}=\{\ev_{h_\b^v}\}_{h_\b^v\in H_\b^v,\,v\in V}$. Then
the EC-setting in the the present case  is established by taking $\mathbf{S}=V$, $\A=P_2\cup H_\b$
and $\fa_{P_2\cup H_\b}=\fa_{P_2}\cup\fa_{H_\b}$. Note that once again  this is an {\it independent}
EC-setting.

\\Clearly all  flaws in $\fa$ are tidy, with seeds of size 1 and with effective size either 1 (if the flaw belongs to $\fa_{P_2}$) or $\b$ (if the flaw belongs to $\fa_{H_\b}$). Thus $E_\A=\{1,\b\}$. As in the previous example $d_1\le \D$.
Moreover,
concerning the estimate of $d_{\b}$, observe that a vertex $w\in V$  has at most $\D$ neighbors
$\{v_1,\dots, v_\D\}$ and therefore $w$ belongs to at most $\D$ distinct neighborhoods $\G^*_G(v_1),\dots \G^*_G(v_\D)$
and in each  $\G^*_G(v_i)$ of these neighborhoods $w$ can be contained  in at most
${\D\choose\b}$
sets of type $h_\b^{v_i}=\{v^i_1,v^i_2.\dots, v^i_{\b+1}\}$. Therefore
\be\label{d1b1}
d_{\b}\le \D {\D\choose\b}\le {\D^{1+\b}\over  \b!}
\ee
Hence we get
$$
\phi_\A=1+\D x+ {\D^{1+\b}\over  \b!}x^\b
$$
\vv
\\The condition (\ref{condentr}) is therefore
\be\label{ex2aa}
k\ge \min_{x>0} {1+\D x +{\D^{1+\b}\over  \b!}x^\b\over x}
\ee
i.e, setting $\D x=y$
\be\label{ex2aa2}
k\ge [\min_{y>0} g_\b(y)]\D
\ee
with
$$
g_\b(y)={1+y+ {\D\over \b!}y^\b\over y}
$$
The minimum occurs in $y_0= \Big[{\D\over \b!}(\b-1)\Big]^{-1\over \b}$ and  $g(y_0)= 1+ \Big({\D\over \b!}\Big)^{1\over \b}(\b-1)^{1\over \b}{\b\over \b-1}$, therefore  there is a $\b$-frugal coloring of  $G$ as soon as
$$
k\ge {\D^{1+{1\over \b}}\over  (\b!)^{1/\b}} {\b\over (\b-1)^{1-{1\over \b}}} +\D
$$
and thus
\be\label{fru}
\ch_\b(G)\le \D+{\D^{1+{1\over \b}}\over  (\b!)^{1/\b}} {\b\over (\b-1)^{1-{1\over \b}}}
\ee

\\The bound (\ref{fru}) represents an improvement
on the previous best bound for the $\b$-frugal chromatic number of a
graph with maximum degree $\D$   given in \cite{NPS}. In
particular, for the case of $2$-frugal coloring (a.k.a. linear
coloring) we get
\be\label{efru}
\ch_2(G)\le { \sqrt{2}}\D^{3\over 2}+\D
\ee
In \cite{Y} it is proved that $\D^{3\over 2}$ is the correct order
since by an explicit example the author proves that $\ch_2(G)\ge
{\D^{3/2}/(6\sqrt{3})}$.

\\To conclude this example let us notice that,  for any $\b\ge 2$ we have once again that $\max_{j=1,\b}{l_j\over j}=2$ and therefore,
for the case  $\b=2$,   the CE criterion (\ref{opa3}) is
\be\label{ex2aa2cl}
k\ge \left[\min_{y>0} {(1+y+ {\D\over 2!}y^2)^2\over y}\right]\D
\ee
yielding
$$
k\ge {8\over\sqrt{27}}\sqrt{2}{\D^{3\over 2}} +{8\D\over 3}+ o(\D)
$$
That is to say in this example as the previous example, EC criterion beats CE criterion by a factor ${8\over\sqrt{27}}\approx 1.54$ in the leading order $\D^{3\over 2}$.

\subsubsection{Acyclic vertex coloring}
A proper coloring $c: V\to [k]$ of the vertices of $G$ is called
an {\it acyclic vertex  coloring} if no cycle of $G$
is bichromatic. The mimumum numbers of colors required
for a graph $G$ to have an acyclic vertex  coloring is called the vertex acyclic chromatic index of $G$ and  it is denoted by $\chi(G)$.

\\Following \cite{AMR},  we say that  a pair of next nearest neighbor vertices $u,v$ of $G$ is a {\it special pair} if
$u$ and $v$ have more than $\a\D^{2/3}$ common neighbors ($\a$ is a positive parameter to be optimized  later).
We set $\mathbf{S}=V$ as the set of atoms and define the set of objects $\A$ as follows. As before $P_2$   is the set of all paths  with $2$ vertices in  $G$; $C_4$ is the set of cycles with four vertices not containing special pairs in  $G$; $S_2$ is the set of all special pairs of vertices in $G$ and finally $P_6$ is the set of all paths with 6 vertices in $G$. We set $\A=P_2\cup C_4\cup S_2\cup P_6$. Let now, for $p\in P_2$, $\ev_p$ be the event that $b$ is monochromatic and let $\fa_{P_2}=\{\ev_p\}_{p\in P_2}$. For $c\in C_4$ let $\ev_c$ be a properly bichromatic 4-cycle  not containing special pairs and let $\fa_{C_4}=\{\ev_c\}_{c\in C_4}$. For $s\in S_2$, let $\ev_s$ be the event that $s$ is monochromatic and let $\fa_{S_2}=\{\ev_s\}_{s\in S_2}$. Finally, for $q\in P_6$, let $\ev_q$ be the event that $p$ is properly bichromatic and let $\fa_{P_6}=\{\ev_q\}_{p\in P_6}$. We set $\fa= \fa_{P_2}\cup \fa_{C_4}\cup \fa_{S_2} \cup\fa_{P_6}$. Clearly a vertex coloring avoiding all
flaws in $\fa$ is acyclic and proper.

\\Note that all events of $\fa$ are tidy and  $\ev\in \fa_{P_2}$ is  such that
$\|\ev\|=1$, $\ev'\in \fa_{C_4}$ is such that $\|\ev'\|=2$, $\ev''\in \fa_{S_2}$ is such that $\|\ev''\|=1$ and $\ev\in \fa_{P_6}$
is such that $\|\ev\|=4$.
Therefore in the present case  we have $E_\A= \{1,2,4\}$.
Moreover let $\tilde d_1$ be the number of paths of size 2 containing a fixed vertex, $\bar d_1$ the set
of special pairs containing a fixed vertex, so that $d_1=\tilde d_1+\bar d_1$; let
$d_2$ be the set of 4-cycles not containing special pairs containing a fixed vertex and let $d_4$
be the set of 6-paths containing a fixed vertex, it is easy to check (see e.g. Lemma 2.4 in \cite{AMR}
or Sec. 2 in \cite{SV})  that $\tilde d_1\le \D$, $\bar d_1\le {1\over \a}\D^{4\over 3}$, $d_2\le {\a\over 2}\D^{8\over 3}$ and $d_4\le 3\D^5$.
Then we have that
$$
\phi_\A(x)= 1+ (\D  +{1\over \a} \D^{4\over 3})x+ {\a\over 2} \D^{8\over 3} x^2 + 3\D^5x^4
$$
Hence the EC criterion  (\ref{condentr}) in the actual case is fulfilled if
\be\label{repve3}
\min_{x>0}{1+ \D x + {1\over \a}\D^{4\over 3}x+ {\a\over 2} \D^{8\over 3} x^2 + 3\D^5x^4\over x}\le k
\ee
Inequality (\ref{repve3}) can be rewritten as
$$
{k}\ge \Big[\min_{x>0}(g_\a(x)+g_1(x))\Big]\D^{4\over 3}
$$
where
$$
g_\a(x)={1+ {1\over \a}x+ {\a\over 2}x^2 \over x} ~~~~~~{\rm and}~~~~~g_1(x) =\D^{-{1\over 3}}(1+x+3x^3)
$$
Taking $x= \sqrt{2\over \a}$ (where   the  minimum  of $g_\a(x)$ occurs) we can bound
$$
{k}\ge (\sqrt{2\a}+{1\over \a})\D^{4\over 3}+ \left[1+\sqrt{2\over \a}+3\left({2\over \a}\right)^{3\over 2}\right]\D
$$
Optimizing now in $\a$ in the leading term $\propto \D^{4\over 3}$ by taking $\a=\sqrt[3]{2}$ we get
\be\label{acyc2}
k\ge {3\over \sqrt[3]{2}} \D^{4\over 3}+ (7+\sqrt[3]{2})\D
\ee
and thus $\chi(G)\le (3/ \sqrt[3]{2}) \D^{4\over 3}+ (7+\sqrt[3]{2})\D$.
This  bound is asymptotically slightly better than the one  obtained in \cite{SV}.

\\To compare this bound with the one obtained via CELL,  observe that  once again  $\max_{j=1,2,4}{l_j\over j}=2$
and therefore  the CE criterion (\ref{opa3}) is
\be\label{ex2aa2cl2}
k\ge \left[\min_{x>0} {(1+ {1\over \a}x+ {\a\over 2}x^2)^2\over x}\right]\D^{4\over 3}+ o(\D^{4\over 3})
\ee
which, optimizing for $\a$ yields
$$
k\ge {9\over 2}{\D^{4\over 3}} + o(\D^{4\over 3})
$$
So in this example  EC beats CELL by a factor $1.89$ in the leading order $\D^{4\over 3}$.

\subsubsection{A dependent EC-setting example: acyclic edge coloring}\label{subacy}
A proper edge coloring of $G$ is \textit{acyclic} if any cycle of $G$ is colored with at least three colors.
The minimum number of colors required such that  $G$ has at least one acyclic proper edge coloring is called
the {\it acyclic edge chromatic number} of $G$  and will be denoted by $a'(G)$.

\noindent If $G$ has maximum degree $\D$, then, for any fixed order  of the edges of $G$, one can properly colors the edges of $E$
using $2\D-1$ colors: just color consecutively
each edge $e$ by using the smallest color not already used in the set $\G(e)$ of the edges adjacent to $e$. Since
$|\G(e)|\le 2\D-2$, at each step this color will be always available. Esperet and Parreau \cite{EP} noticed
that, still using just $2\D-1$ colors, one can also avoid  at each step
all bichromatic cycles of length 4. Indeed, when coloring $e$, if $f$ and $f'$ are in $\G(e)$ and have the same
color $s_1$ and $e'\in E$, having color $s_2$, is  such that
$\{e,f,f',e'\}$ is a 4-cycle, then we must avoid the color $s_2$ but at the same time in $\G(e)$ only at most $2\D-3$ colors
are used (since $f$ and $f'$  have the same color). In other words the following lemma holds.

\begin{lem}\label{4cy} Given a graph $G=(V,E)$ with maximum degree $\D$, the edges of $G$ can be colored sequentially
using just $2\D-1$ colors, in such a way that
the coloring is proper and  contains no  bichromatic cycle of length 4.
\end{lem}

\\In order  to use the EC criterion to  estimate $a'(G)$, let us first of all establish the pair $(\mathbf{S},k)$.
Since here we are coloring the edges of $G$ we set $\mathbf{S}=E$ and, in view of Lemma \ref{4cy},we choose  $k$ in such a way that
$k\ge 2\D-1$.

\\Now, for any partial coloring $\g\in [k]_0^E$ of the edges of $G$ using $k\ge 2\D-1$ colors
and for any $e\in E$, let $\F(e,\g)$ be the subset of $[k]$ formed by the colors  present in $\G(e)$  or such that
exists $f,f'\in \G(e)$  and $e'\in E$ such that
$\{e,f,f',e'\}$ is a properly bichromatic 4-cycle in $G$.  By the discussion above we have that $|\F(e,\g)|\le 2(\D-1)$.
We now  associate  to each $e\in E$  and each $s\in  \F(e,\g)$  the {\it non-standard} elementary flaw $\ev=\{\mbox{$e$ has the color $s$}\}$.
Let us denote by $\fa_E$ the set of
all these flaws. We stress once again that each flaw  $\ev\in\fa_E$ is non-standard. Moreover, since all events
of $\fa_E$ are elementary, we have  $\|\ev\|=|\obj(\ev)|=1$ for all $\ev\in\fa_E$ . Clearly a coloring $c: E\to [k]$  avoiding all   flaws in $\fa_E$
is a proper coloring of the edges of $G$ with no bichromatic 4-cycles.
Let now $C_{2n}$ be the set of all cycles of size $2n$ in $G$ and let $C=\{C_{2n}\}_{n\ge 3}$ the set of all cycles of even length greater than 4 in $G$. Let $\fa_C$
the family of flaws  $\fa_C=\{\mbox{the cycle $c$ is properly bi-chromatic}\}_{c\in C}$. Any flaw of $\ev\in \fa_C$ is standard,
non elementary and  tidy with seeds of size 2. Thus, for any $\ev\in \fa_C$ and
 for any $e\in \obj(\ev)$,  it holds  $\|\ev_e\|=|\obj(\ev)|-{2}$.

\\Clearly a coloring $c: E\to [k]$  avoiding all   flaws in $\fa_E\cup \fa_C$
is an acyclic edge coloring of $G$.
Therefore the EC-setting in the the present case  is established by taking $\A=E\cup C$
and $\fa=\fa_E\cup\fa_C$. We stress once again that this is a {\it dependent} EC-setting since each  event $\ev$ in the subfamily  $\fa_E$ is non standard (i.e.
$\ev$
depend on the (partial) coloring outside $\obj(\ev)$).
Moreover we have
$$
E_\A= \{1\}\cup \{2n-2\}_{n\ge 3}
$$
Let us find the numbers $d_j$ defined in (\ref{dj}) for $j\in E_\A$.
Recalling that
$$
\max_{v\in V,\g\in [k]_0^V}|\F(v,\g)|\le 2(\D-1)
$$
we get
$$
d_1\le 2(\D-1)
$$
Secondly, for $n\ge 3$,
the number $d_{2n-2}$  represents here
the maximal number of cycles  of size $2n$ in $G$ containing a fixed edge and it can be bounded as
$$
d_{2n-2} \leq (\D-1)^{2n-2},
$$
so that, for $0<x<\D-1$, we have
$$
\phi_\A(x)\le 1+2(\D-1) x+\sum_{n\ge 3}(\D-1)^{2n-2}x^{2n-2} = 1+2(\D-1) x+ {[(\D-1)x]^4\over 1- [(\D-1)x]^2}
$$
Hence the entropy-compression criterion  (\ref{condentr}) in the actual case is fulfilled if
\be\label{repve4}
k\ge \min_{x>0}{ 1+2(\D-1) x+ {[(\D-1)x]^4\over 1- [(\D-1)x]^2}\over x}
\ee
with $x>0$.
Setting $(\D-1) x=y$, inequality (\ref{repve4}) can be rewritten as
$$
{k}\ge \left(\Big[\min_{y>0}g(y)\Big]\right)(\D-1)
$$
where
$$
g(y)={1+2y+ {y^4\over 1- y^2}\over y}
$$
The $\min_{y>0}g(y)$ is attained at $x= {1\over 2}(\sqrt{5}-1)$ and $g({1\over 2}(\sqrt{5}-1))=4$
and thus
\be\label{astar2}
a'(G)< 4(\D-1)
\ee
which is the bound given by Esperet and Parreau in   \cite{EP}.

\vskip.2cm
\\{\bf Remark}. As said in the introduction, the constant  4 here above is sensibly better than 9.163
obtained via the cluster expansion criterion in \cite{NPS}. However this strong improvment is much
more due to Lemma \ref{4cy} that to entropy compression. To see this, first observe the cluster
expansion criterion applied to the acyclic edge coloring in the ways illustrated  in \cite{NPS} is
obtained by  assigning to each edge of $G$ a color chosen  uniformly and independently between  $N$ available colors.
In view of Lemma \ref{4cy} we can adopt a more  efficient strategy. Suppose that our random coloring of the edges of $G$
is obtained
by using $N=2(\D-1) +k$ colors (with $k\ge 1$) and by coloring sequentially the  edges of $G$,
in such a way that the edge $e$ colored at step $t$ is  colored by choosing uniformly and independently a number $s\in [k]$
and then coloring $e$ using the $s^{th}$ available color among those avoiding monochromatic cherries  and bichromatic 4-cycles (which, at each step, are at least $k$). With this kind of random coloring the possible bad events  are only those in the family $\fa_C$, i.e. they are the
bichromatic cycles
of size greater or equal than 6. Due to our random coloring we  are not in the Moser-Tardos setting but nothing prevent us to  apply the general Theorem \ref{bis}. Let us choose as graph $\GG$ the graph with vertex set $\fa_C$ and
edge set $\E$ such that $\{\ev,\ev'\}\in \E$ if and only if $\obj(\ev)\cap\obj(\ev')\neq 0$. It is not difficult to check that,
given the random greedy coloring of the edges of $E$ decribed above, for any $\ev\in \fa_C$ and any $U\in \fa\setminus\G_\GG(\ev)$, we have
$$
\prob\left(\ev\Big|\bigcap_{\ev'\in U}\bar{\ev'}\right)\;\leq\;{1\over k^{2|\obj(\ev)|-2}}
$$
Therefore, using bounds (\ref{des})-(\ref{opa3}), condition (\ref{cond1n})  can be written as
\be\label{opacy}
k~\ge~~ \min_{0<x<(\D-1)^{-1}}{\left[1+\sum_{k\ge 3}(\D-1)^{2k-2}x^{2k-2}\right]^{3\over 2}\over x}
\ee
where $3/2=\sup_{\ev\in \fa_C}|\obj(\ev)|/(|\obj(\ev)|-2)$.
Condition (\ref{opacy}) is satifsied as soon as
$$
k\ge\left[ \min_{0<y<1}{\left(1+ {y^4\over (1-y^2)}\right)^{3\over 2}\over y}\right](\D-1)= 2.181(\D-1)
$$
Therefore  it is possible to color the edges of $G$ using $N\ge 2(\D-1)+ 2.181(\D-1)=4.181(\D-1)$ in such a way that
the resulting coloring is proper and with no bichromatic cycles of any size, This bound  is very close to that obtained
by Esperet and Parreau. In conclusion, bearing in mind  Lemma \ref{4cy},  as far as acyclic edge coloring is considered,
CELL is beaten by EC only because of the exponent $3/2$ in the numerator of the r.h.s. of (\ref{opacy}).

\\We finally want to stress that Giotis et al. \cite{GKPT} have recently improved the Esperet-Parreau bound proving that $a'(G)\le 3.79(\D-1)$
combining Lemma \ref{4cy} with a variant of the Moser-Tardos algorithm specifically designed for the case of the acyclic edge coloring.

\subsection{EC and CELL can tie the game. A non graph coloring example: Independent sets}
\\Let $G$ be a graph with vertex set $V$, edge set $E$ and with maximum degree $\D$.
We want to find the least integer  $k_\D$ such that, for any partition $\{V_1,\dots, V_n\}$ of
$V$ with $|V_i|\geq k$,  there exists an independent subset $I$ of $V$  with
exactly one vertex in each $V_i$. This set $I$ is called  an  {\it independent tranversal} of $G$ w.r.t. the partition $\{V_1,\dots, V_n\}$.
It is long known that $k_\D$ is linear in  $\D$ \cite{Al,Fe}. The best bound for $k_\D$, i.e. $k_\D\le 2\D$ is due
to Penny Haxell \cite{H}. The Haxell's bound is tight  since Szab\'o and Tardos \cite{ST} showed that $k_\D> 2\D-1$.

\\Let us see in this section  what EC and LLL are able to say about this problem.
Without loss of generality, we can assume that $G$ is such that $|V_i| = k$ for all
sets $V_i$. The general case  $|V_i|\geq k$ follows
by considering the graph induced by $G$ on a union of $n$ subsets
of cardinality $k$, each of them a subset of one $V_i$.

\\The EC setting in this case can be established as follows.
We may assume that  for each $i\in [n]$ the set $V_i$ is ordered,
namely, $V_i=\{v_1^{(i)}, \dots , v_k^{(i)}\}$. If now we identify, for each $j\in [k]$ and
$i\in [n]$,  $v_j^{(i)}$ with $j$,   we can take as  set of atoms $\mathbf{S}=\{V_1,V_2,\dots, V_n\}$, and the
 set of ``colors"  can be identified  with $[k]$,  a ``coloring" $c$ of $\mathbf{S}$ being just the
 selection of an element in each $V_i$. Namely, to choose the color $s\in [k]$ for the atom $V_i$ means
 to select the vertex $v^{(i)}_s\in V_i$.
Given $i,j\in [n]$, we say that the pair $\{V_i,V_j\}$ is connected in $G$ if there exists
$v\in V_i$ and $v'\in V_j$ such that $\{v,v'\}\in E$.
Then the family of objects $\A$ in the present case
is constituted by the connected pairs $\{V_i, V_{i'}\}$ in $G$. Given $A=\{V_i, V_{i'}\}\in \A$ let
$\ev$ be the event in $[k]^A$
$$
\ev=\{c(V_i)=v,\; c(V_j)=v' \;\mbox{and} \;\{v,v'\}\in E\}
$$
Note that $\ev$ is standard and elementary.
Let now, for each $A=\{V_i, V_{i'}\}\in \A$,  $\fa_A$  be the family of flaws
$$
\fa_A=\Big\{\ev: c\in [k]^A:\;\{c(V_i), c(V_j)\}\in E\Big\}
$$
In other words, for each $A\in \A$, $\fa_A$ is constituted  by  the (elementary) events formed by the colorings that select a vertex $v$ in $V_i$
and a  vertex $v'$ in $V_{i'}$ and $\{v,v'\}$ is an edge of $G$. The family
$\fa=\cup_{A\in \A}\fa_A$,  constituted by standard elementary flaws, is such that a
coloring of $\mathbf{S}$ avoiding all the events in the family $\fa$ selects an independent set of $G$.

\\In the present case, since
any $\ev\in \fa_A$ is elementary, we have that $\|\ev\|=|\obj(\ev)|=2$ for all $\ev\in \fa$.
This means that $E_\A=\{2\}$.
Finally  we have to estimate  $d_2$, i.e., the maximal number of events  containing a fixed $V_i\in \A$.
Since any $V_i$ has $k$ vertices
and each of these vertices is adjacent to at most $\D$ other vertices of $G$, we get that
$d_2\leq k\D$ and thus
$$
\phi_\EA(x)=1+ k\D x^2
$$
Therefore condition (\ref{condentr}) is written as follows
$$
k \ge 2(k\D)^{\frac{1}{2}}
$$
which is to say $k_\D\le  4\D$. This bound is  equal to that obtained with the CELL (see e.g. \cite{BFPS}). By the discussion
made in Section \ref{compar} this is not a surprise. Indeed, in the present case
EC criterion and CELL criterion yield the same bound because any $\ev\in \fa$ is elementary.

\\Finally it is worth to mention that very recently  Graf and Haxell \cite{GH} have designed an efficient algorithm (not related to EC or Moser-Tardos algorithmic LLL!) able to find an independent  transversal in a graph $G=(V,E)$ with maximum degree $\D$ for any partition
$\{V_1,\dots, V_n\}$ of $V$ such that $|V_i|\ge 2\D+1$ for $1=1,\dots, n$.

\subsection*{Declarations of interest: none}

\section*{Acknowledgments}
A.P.  and R.S. have been partially supported by the Brazilian  agencies
Conselho Nacional de Desenvolvimento Cient\'\i fico e Tecnol\'ogico
(CNPq) and  Funda\c c\~ao de Amparo \`a  Pesquisa do Estado de Minas Gerais (FAPEMIG - Programa de Pesquisador Mineiro).


\begin{thebibliography}{99}

%\bibitem{AG} Achlioptas, D. Gouleakis, T.: {\it Algorithmic Improvements of the Lov\'asz Local Lemma via Cluster Expansion}.
%Annual Conference on Foundations of Software Technology and Theoretical Computer Science, FSTTCS, 16-23, (2012).


%\bibitem{AIV} Achlioptas, D.; Iliopoulos, F.; Vlassis, N.: {\it Stochastic Control via Entropy Compression},
%Preprint, arxiv.org/abs/1607.06494 (2016).


\bibitem{Al} Alon, N.: {\it The linear arboricity of graphs}, Israel Journal of Mathematics, {\bf 62}, 311-325 (1988).

\bibitem{A} Alon, N.: {\it A parallel algorithmic version of the local lemma}. Random Structures and Algorithms,
{\bf 2}, n. 4, 367-378  (1991).

\bibitem{AMR} Alon, N. Mc Diarmid, C.; Reed, B.: {\it  Acyclic colouring of graphs}, Random Structures and Algorithms, {\bf 2}, no. 3, 277-288 (1991).

\bibitem{AS} Alon, N. and Spencer, J.: {\it The Probabilistic Method. Third  Edition}. New York, Wiley-Interscience, (2008).

%\bibitem{AP} Alves, R. G.; Procacci, A.:{\it Witness trees in the Moser-Tardos algorithmic Lov\'asz Local Lemma and
%Penrose trees in the hard-core lattice gas}, Journal of Statistical Physics, {\bf 156}, p. 877-895 (2014)

\bibitem{B} Beck, J.: {\it An Algorithmic Approach to the Lov\'asz Local Lemma}, Random Structures and Algorithms,
{\bf 2}, n. 4, 343-365 (1991).


\bibitem{Be} Bernshteyn, A.:  {\it The local cut lemma},  European Journal of Combinatorics, {\bf 63}, 95-114 (2017).


\bibitem{BFPS} Bissacot, R.; Fern\'andez, R.; Procacci A.; Scoppola, B.: {\it An Improvement of the Lov\'asz Local Lemma via Cluster Expansion},
Combinatorics Probability and Computing, {\bf 20}, n. 5, 709-719
(2011)

\bibitem{BKP}  B\"{o}ttcher, J.; Kohayakawa, Y.; Procacci, A.: {\it Properly coloured copies and rainbow
copies of large graphs with small maximum degree},
Random Structures and Algorithms, {\bf 40}, n. 4,  425-436 (2012).

\bibitem{CLY} Cai, J.; Li, X.; Yan, G.: {\it Improved upper bound for the degenerate and star
chromatic numbers of graphs} J. Comb. Optim.,  {\bf 34}, 441-452 (2017).

\bibitem{CR} Camungol, S.; Rampersad, N.:  {\it Avoiding approximate repetitions with respect to the
longest common subsequence distance}, Involve, a Journal of Mathematics, {\bf 9},  No. 4, 657-666 (2016).

\bibitem{DJKW} Dujmovic, V.; Joret, G.; Kozik, J.; Wood, D. R.:
{\it Nonrepetitive colouring via entropy compression}, Combinatorica, {\bf 36}, n.6,  661-686 (2016).

\bibitem{Dm} Drmota, M.: {\it Combinatorics and asymptotics on trees}, Cubo J. {\bf 6} (2) (2004).

\bibitem{D} Dobrushin , R. L.: {\it Perturbation methods of the theory of
Gibbsian fields}, in P. Bernard (editor), Lectures on Probability Theory and Statistics,
Lecture Notes in Mathematics Volume 1648, 1996, pp 1-66
Springer-Verlag, Berlin, (1996).

\bibitem{EL} Erd\H{o}s, P. and Lov\'{a}sz, L.: {\it Problems and results on 3-chromatic
hypergraphs and some related questions, in Infinite and finite
sets}. Vol. II, Colloq. Math. Soc. Janos Bolyai, Vol. 10, pp. 609-627.
North-Holland, Amsterdam, (1975).

\bibitem{EP} Esperet, L.; Parreau, A.: {\it Acyclic edge-coloring using entropy
compression}, European Journal of Combinatorics, {\bf 34}, In. 6,
1019- 1027, (2013).

\bibitem{Fe} Fellows, M.: {\it Transversals of vertex partitions in graphs}, SIAM Journal of Discrete Mathematics
{\bf 3}, 206-215, (1990).

\bibitem{FP} Fern\'andez, R.; Procacci A.: {\it Cluster
expansion for abstract polymer models. New bounds from an old
approach}, Communications in Mathematical Physics. {\bf 274}, n.1,
123-140 (2007).

\bibitem{FS} Franceti\'c, N.;  Stevens, B.:
{\it Asymptotic size of covering arrays: an application of entropy compression},
J. Comb. Designs, {\bf  25}, Issue 6,
243-257 (2017).

%\bibitem{GH}  Gasarch W.; Haeupler B.: {\it Lower Bounds on van der Waerden Numbers: Randomized and
%Deterministic-Constructive}, The Electronic J. Combinatorics, {\bf 18} (2011), $\#$P64.

\bibitem{GKPT} Giotis, I.; Kirousis, L.; Psaromiligkos, K. I.; Thilikos, D. M.: {\it  Acyclic edge coloring
through the Lov\'asz Local Lemma}, Theoretical Computer Science, Elsevier, 665, pp.40 - 50, (2017).

\bibitem{GH} Graf, A.; Haxell. P.: {\it Finding Independent Transversals Efficiently}, arXiv:1811.02687 (2018).

\bibitem{GKM} Grytczuk, J.; Kozik, J.; Micek, P.: {\it New approach to nonrepetitive sequences},
Random Struct. Algorithms, {\bf 42}, Issue 2, 214-225 (2013).

\bibitem{GMP} Gon\c{c}alves, D.; Montassier, M.; Pinlou, A. {\it Entropy compression
method applied to graph colorings}, arXiv:1406.4380 (2014).


%\bibitem{HS} Harris, D. G.; 	Srinivasan, A. : {\it Algorithmic and Enumerative Aspects
%of the Moser-Tardos Distribution}, {\bf 13}, Issue 3,
%Article No. 33 (2017).

%\bibitem{HSS} 	Haeupler, B.;
%Saha, B.; 	Srinivasan, A.: 	{\it New Constructive Aspects of the Lov\'asz Local Lemma}, {\bf 58}, Issue 6,
%Article No. 28 (2011).

\bibitem{H} Haxell, P. E.: {\it  A note on vertex list colouring}. Combin. Probab. Comput. {\bf 10}, no. 4, 345-347 (2001).

%\bibitem{KS} Kolipaka, K. B. R.; Szegedy, M., {\it Moser and Tardos meet Lov\'asz},
%Proceedings of the 43rd annual ACM symposium on Theory of computing
%Pages 235-244, ACM New York, NY, USA (2011).

%\bibitem{MD} McDiarmid. C.: {\it Hypergraph coloring and the Lov\' asz Local Lemma},
%Discrete Math., {\bf 167/168}, 481-486  (1997).

\bibitem{Mo} Moser, R. A.: {\it A constructive proof of the Lov\'asz Local Lemma}, in Proceedings of the 41st Annual
ACM Symposium on the Theory of Computing (STOC). ACM, New York (2009).

\bibitem{Mo2} Moser, R. A.: {\it Exact algorithms for constraint satisfaction problems},
Diss. ETH Z\"urich, Nr. 20668, Logos Verlag, Berlin (2013).

\bibitem{MT} { Moser, R. ; Tardos, G.: {\it A constructive proof of the general Lov\'asz Local
Lemma}, J. ACM {\bf 57} Article 11, 15 pages  (2010).}

\bibitem{NPS} Ndreca, S.;  Procacci, A.; Scoppola, B.: {\it
Improved bounds on coloring of graphs},  European Journal of
Combinatorics, {\bf 33}, n 4, p. 592-609 (2012).


\bibitem{OP} Ochem, P.; Pinlou, A.: {\it Application of Entropy Compression in Pattern Avoidance},
Electr.  J. Comb. {\bf  21}, Issue 2 (2014).

\bibitem{Pe} Pegden, W.: {\it An extension of the Moser-Tardos algorithmic local lemma},
SIAM J. Discrete Math. {\bf 28}, 911-917 (2013).

\bibitem{PS} Procacci, A.; Sanchis, R.: {\it Perfect and separating hash families:
new bounds via the algorithmic cluster expansion local lemma},
Annales  de l'Institut Henry Poincar\'e D  combinatorics,
physics and their interactions, to appear (2016).

\bibitem{Pr} Przybylo, J.: {\it On the Facial Thue Choice Index via Entropy Compression},
Journal of Graph Theory, {\bf 77}, Issue 3, 180-189, (2014).

\bibitem{PSS} Przybylo, J.; Schreyer, J.; $\check{S}$krabul'\'akov\'a, E.: {\it On the facial Thue number
of plane graphs Index via Entropy Compression}, Graphs and Comb. {\bf 32}, n.3, 1137-1153 (2016).

%\bibitem{PS} Procacci A.; Sanchis R.: {\it Perfect and separating
%Hash families: new bounds via the algorithmic cluster expansion local lemma}. Preprint arXiv:1601.05389.

%\bibitem{RS} Radhakrishnan, J. ;  Srinivasan, A.; {\it Improved bounds and algorithms
%for hypergraph two-coloring}, Random Struct. Algorithms
%{\bf 16}, 4-32, (2000).

\bibitem{SS}  Scott, A.; Sokal, A. D.: {\it The repulsive lattice gas,
the independent-set polynomial, and the Lov\'asz local lemma},
 J. Stat. Phys.  {\bf 118}, no. 5-6, 1151--1261, (2005).

 \bibitem{SV} Sereni, J. S.; Volec, J.:{\it  A note on acyclic vertex-colorings}. J. Comb., {\bf 7}(4), 725-737, (2016).

 \bibitem{Sh} Shearer, J. B.: {\it On a problem of Spencer}. Combinatorica {\bf 5}, 241-245, (1985).

 \bibitem{ST} Szab\'o, T.; G. Tardos, G.: {\it Extremal problems for transversals in graphs with bounded degree},
Combinatorica. {\bf 26}(3),  333-351 (2006).

\bibitem{Ta} Tao, T.: {\it Moser's entropy compression argument},  Terence Tao Blog post, (2009).

\bibitem{Y} Yuster, R.: {\it Linear coloring of graphs}, Discr. Math. {\bf 185}, Issues 1-3, Pages 293-297 (1998).




%\bibitem{Sh} Shearer, J. B.: {\it On a problem of Spencer},
%Combinatorica {\bf 5}, 241-245, (1985).

\end{thebibliography}
\end{document}